\documentclass[a4paper,11pt]{article}
 \pdfoutput=1

\usepackage{floatrow}
\usepackage{mathtools}
\usepackage{hyperref}
\usepackage{amssymb}
\usepackage{tikz,pgfplots}
\usepackage{algorithm}
\usepackage[noend]{algpseudocode}
\usepackage{amsfonts,amssymb,amsmath,amsthm}
\usepackage{color,ulem}
\usepackage{fullpage}
\usepackage{graphicx}
\usepackage[active]{srcltx}
\usepackage{mathabx}
\usepackage{verbatim}
\usepackage{cancel}
\usepackage{mathtools}
\usepackage{tikz,pgfplots}
\usepackage{natbib}
\usepackage{tabu}
\usepackage{cite}
\usepackage{hyperref}
\hypersetup{
  linkcolor  = green!60!black,
  citecolor  = blue!60!black,
  urlcolor   = red!60!black,
  colorlinks = true,
}
\usepackage{caption}
\usepackage{subcaption}
\usepackage[official]{eurosym}

\newtheorem{assumption}{Assumption}

\def\argmin{\mathop{\rm argmin}}
\newtheorem*{conjecture*}{Conjecture}

\newtheorem{remark}{Remark}

\newfloatcommand{capbtabbox}{table}[][\FBwidth]
\usepackage[T1]{fontenc} 
\usepackage[utf8]{inputenc}
\usepackage{times}
\usepackage{a4wide}
\usepackage{amsthm, amsmath, amssymb}
\usepackage{algorithm}
\usepackage{float}
\usepackage{subcaption}
\usepackage[noend]{algpseudocode}
\usepackage{cases}
\usepackage{ dsfont }
\usepackage{array}

\usepackage{chngcntr}
\usepackage{mathrsfs}
\usepackage{hyperref}
\usepackage{fancyvrb}
\usepackage{bbm}
\usepackage{bm}
\usepackage{latexsym}
\usepackage{multicol}
\usepackage{eucal}
\usepackage[font={small}]{caption}
\usepackage{lipsum}                    	
\usepackage{verbatim}
\usepackage{graphicx}
\usepackage{grffile}
\usepackage{caption}
\usepackage{subcaption} 
\usepackage{url}
\usepackage{authblk}

\def\argmin{\mathop{\rm argmin}}

\newcommand{\CK}{\mathcal{K}}

\definecolor{abblue}{rgb}{0,0.2,0.6}

\definecolor{abgreen}{rgb}{0,0.4,0.1}

\begin{document}

\title{Regression Monte Carlo for Microgrid Management}

\author[1]{Clemence Alasseur\thanks{email: clemence.alasseur@edf.fr}}\affil[1]{EDF R\&D  - FIME, Palaiseau, France; }
\author[2]{Alessandro Balata\thanks{email: A.Balata@leeds.ac.uk}}\affil[2]{University of Leeds, Woodhouse Lane, Leeds LS2 9JT, United Kingdom; }
\author[3]{Sahar Ben Aziza\thanks{email: sahar.benaziza@enit.utm.tn}}\affil[3]{University of Tunis El Manar,ENIT-LAMSIN, BP.37, Le Belvédère 1002 Tunis, Tunisia;}
\author[4]{Aditya Maheshwari\thanks{email: aditya\textunderscore maheshwari@umail.ucsb.edu}}\affil[4]{University of California, Santa Barbara, USA;}
\author[5]{Peter Tankov\thanks{email: peter.tankov@ensae.fr}}\affil[5]{CREST-ENSAE, Palaiseau, France;}
\author[1]{Xavier Warin\thanks{email: xavier.warin@edf.fr}}
\maketitle

\begin{abstract} 
We study an islanded microgrid system designed to supply a small
village with the power produced by photovoltaic panels, wind turbines
and a diesel generator. A battery storage system device is used to
shift power from times of high renewable production to times of high
demand. We build on the mathematical model introduced in
\citet{heymann17} and optimize the diesel consumption under a
``no-blackout'' constraint. We introduce a methodology to solve
microgrid management problem using different variants of Regression
Monte Carlo algorithms and use numerical simulations to infer results about the optimal design of the grid. 
\end{abstract}

\section{Introduction}
A Microgrid is a network of loads and energy generating units that
often include renewable sources like photovoltaic (PV) panels and wind
turbines alongside more traditional forms of thermal electricity
production. These microgrids can be part of the main grid or
isolated. Communities in
rural areas of the world have long now enjoyed the installation of
isolated microgrid systems that provide a reliable and often environment-friendly
source of electricity to meet their power  needs. 

The elementary purpose of a microgrid is to provide a continuous
electricity supply from the variable power produced by renewable
generators while minimizing the installation and running costs. In
this kind of systems, the uncertainty of both, the load and the
renewable production is high and its negative effect on the system
stability can be mitigated by including a  battery energy storage system
in the microgrid. Energy storage devices ensure power quality,
including frequency and voltage regulation (see \citet{hayashi2017})
and provide backup power in case of any  contingency. A dispatchable unit in the form of diesel generator is also used as a backup solution and to provide baseload power. 

In this paper, we consider a traditional microgrid serving a small group of customers in islanded mode, meaning that the
network is not connected to the main national grid.
The system consists of an
intermittent renewable generator unit, a conventional dispatchable
generator, and a battery storage system. Both the load and the
intermittent renewable production are stochastic, and we use a
stochastic differential equation (SDE) to model directly the residual
demand, that is, the difference between the load and the renewable
production. We then set up a stochastic optimization problem, whose
goal is to minimize the cost of using the diesel generator plus the
cost of curtailing renewable energy in case of excess production,
subject to the constraint of ensuring reliable energy supply. A
regression Monte Carlo method from the mathematical finance literature
is used to solve this
stochastic optimization problem numerically. Three variants of the regression alrogithm, called grid
discretization, Regress now and Regress later are proposed and
compared in this
paper. The numerical examples illustrate the performance of the
optimal policies, provide insights on the optimal sizing of the
battery, and compare the policies obtained by stochastic optimization
to the industry standard, which uses deterministic policies. 

The optimization problem arising from the search for a cost-effective
control strategy has been extensively studied. Three recent survey papers \citet{olivares2014trends,Reddy2017,Liang2014} summarize different methods used for optimal usage, expansion and voltage control for the microgrids.  Heymann
et. al.\citet{heymann16,heymann17}  transform the optimization problem
associated with the microgrid management into an optimal control
framework and solve it using the corresponding Hamilton Jacobi Bellman
equation. Besides proposing an optimal strategy, the authors also
compare the solution of the deterministic and stochastic
representation of the problem. However, similarly to most PDE methods,
this approach suffers from the curse of dimensionality and as a
result, it is difficult to scale. The main contribution of this paper
is to solve the microgrid control problem using Regression Monte Carlo
algorithms. In contrast to existing approaches, the method
used in this paper is more easily scalable and works well in
moderately large dimensions \citet{bouchard2012monte}. 
%In this paper, we convert the
%optimal control problem into a statistical learning problem using
%dynamic programming principle. This transformation allows us to use the Regression Monte Carlo methods from the mathematical finance literature to solve the microgrid management problem presented here.

Identifying the optimal mix, the size and the placement of different
components in the microgrid is an important challenge to its large
scale use. The papers \citet{Mashayekh2017a, Mashayekh2017b} use mixed-integer
linear programming to address the design problem and test their model
on a real data set from a microgrid in Alaska. In a similar work,
\citet{Olatomiwa2015} studied the economically optimal mix of PV, wind,
batteries and diesel for rural areas in Nigeria. In
\citet{haessig2015energy}, optimal battery storage sizing is deduced
from the autocorrelation structure of renewable production forecast
errors. 
In this paper, we
propose an alternative approach for the optimal sizing of the battery
energy storage system, assuming stochastic load dynamics and fixed
lifetime of the battery. Our in-depth analysis of the system behavior leads to practical guidelines for the design and control of islanded microgrids. 

Finally, several authors
\citet{ding2012stochastic,ding2015rolling,collet2017optimal} used
stochastic control techniques to determine optimal operation strategies for wind
production -- storage systems with access to energy markets. In
contract to these papers, in the present study, energy prices appear
only as constant penalty factors in the cost functional, and the main focus
is on the stable operation of the microgrid without blackouts. 

The rest of the paper is organized as follows:
In section \ref{model_description} we describe the microgrid model and introduce the different components of the system, in section \ref{control_problem} we translate the problem of managing the microgrid in a stochastic optimization problem and present the dynamic programming equation that we intend to solve numerically. Section \ref{numerical_methods} introduces the numerical algorithms used to solve the control problem, we give a general framework for solving the dynamic programming equation and we then provide three algorithms for the approximation of conditional expectations. 
In section \ref{experiments} we illustrate the results of the numerical experiments, identify the best algorithm among those we studied and then employ it to analyze the system behavior.
We conclude with section \ref{deterministic_comparison} where the estimated policy for the stochastic problem is compared, in an appropriate manner, with a deterministically trained one; the aim is to provide evidence that industry-widespread deterministic approaches underperform stochastic methods.

\section{Model description}
\label{model_description}

In this section, we will discuss the topology of the microgrid, its
operation, components and their respective dynamics. Although we
discuss a simplified microgrid model, more complicated typologies can be studied using straightforward generalizations of the methods presented in this paper.

Consider a microgrid serving a small, isolated village; most of the power to the village is supplied by generating units whose output has zero marginal cost, is intermittent and uncontrolled.  Additional power is supplied by a controlled generator whose operations come alongside a cost for the microgrid owner (either the community itself or a power utility). Often the intermittent units include PV panels and wind turbines, while the controlled unit is often a diesel generator. In order to fully exploit the free power generated by the renewable units at times when production exceeds the demand, microgrids are equipped with energy storage devices. These can be represented by a battery energy storage system.

The introduction of the battery in the system not only allows for inter-temporal transfer of energy from times when demand is low, to times when it is higher, but also introduces an element of strategic behavior that can be employed by the system controller, to minimize the operational costs. Without an energy storage,  diesel had to be run at all times demand exceeded production. When a battery is installed, intensity and timing of output from the diesel generator can be adjusted to move the level of charge of the battery towards the most cost effective levels.

\begin{figure}
\includegraphics[width=0.6\textwidth ]{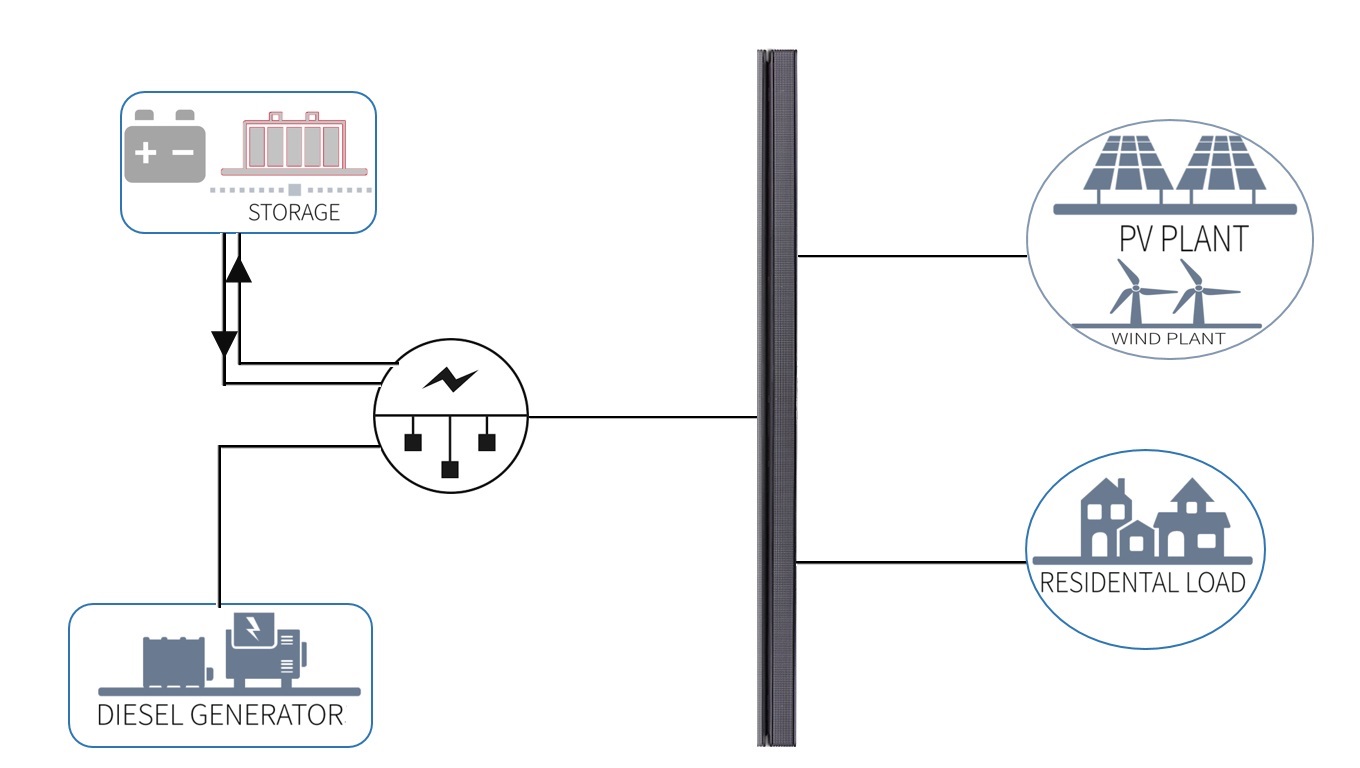}
\caption{The figure above shows an example of microgrid topology that contains all the elements in our model. The network is arranged as follows: photovoltaic panels and wind turbines provide renewable generation, a diesel generator provides dispatchable power for the village and a battery storage system is used to inject or withdraw energy. }
\label{fig:microgrid_topology}
\end{figure}

In figure \ref{fig:microgrid_topology} we propose a schematic description of the system which might help the reader to familiarize themselves with the microgrid,  whose components are described more in depth in the following subsections.

\begin{remark}
Note that for convenience, in the following, we will work in discrete time only. This setting is not restrictive as in reality measurements of the systems are repeated at a given, finite, frequency.
We also consider a finite optimization horizon represented by the number of periods over which we want to optimize the system operations indicated by $T$
\end{remark}

\subsection{Residual Demand}

Consider two stochastic processes $L_t$ and $R_t$, the former represents the demand/load and the latter the production through the renewable generators. Notice that both processes are uncontrolled and they represent, respectively, the unconditional withdrawal or injection of power in the system (constant during time step). For the purpose of managing the microgrid, the controller is interested only in the net effect of the two processes denoted by the process $X_t$: 
\begin{equation}
\label{eq:residual_demand}
X_t=L_t-R_t\; ;\quad t\in\{0,\,1,\,\dots,\,T\}.
\end{equation}

\begin{remark}
The state variable $X_t$ represents the residual demand of power at each time $t$, such that for $X_t>0$, we should provide power through the battery or diesel generator and for $X_t<0$  we can store the extra power in the battery.
\end{remark}

For simplicity, we model the residual demand as an AR(1) process, the discrete equivalent of an Ornstein–Uhlenbeck process. In practical applications we expect $X_t$ to be an $\mathbb{R}$-valued mean reverting process with many different sources of noise and time dependent random parameters; our formulation avoids the cumbersome notation using constants in place of stochastic processes still providing scope for generalization.
The process $X_t$ is driven by the following difference equation, starting from an initial point $X_0=x_0$: 
\begin{equation}
\label{eq:residual_demand_dynamics}
X_{t+1}= X_{t}+ b(\Lambda_{t} -X_{t} )\Delta t + \sigma \sqrt{\Delta t}\; \xi _{t} \; ; \quad t\in\{0,1,\dots,T\}
\end{equation}

where $\xi_{t}\sim\mathcal{N}(0,1)$, $\Delta t$ is the amount of time before new information is acquired, $b$ is the mean reversion speed, $\sigma$ the volatility of the process and $\Lambda_t$ is the time dependent mean reversion level.

\begin{remark}In real applications the function $\Lambda_t$ should represent the best forecast available for future residual demand at the time of the estimation of the policy. 
\end{remark}

\subsection{Diesel generator}

The Diesel generator represents the controlled dispatchable unit. The state of the generator is represented by $m_t=\{0,1\}$. If $m_t=0$ then the diesel generator is OFF, while it is ON when $m_t=1$. When the engine is ON, it produces a power output denoted by $d_t \in [d_{min},d_{max}]$ at time $t$, for $d_{min}>0$.

\begin{figure}[H]
        \centering
        \begin{subfigure}[b]{0.45\textwidth}
            \centering
            \includegraphics[width=\textwidth]{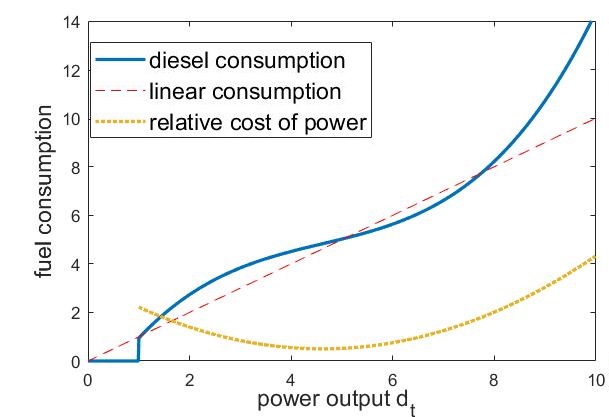}
            \caption{$\rho(d)=\frac{(d-6)^3+6^3+d}{10}$}    
        \end{subfigure}
        \quad
        \begin{subfigure}[b]{0.45\textwidth}  
            \centering 
            \includegraphics[width=0.95\textwidth]{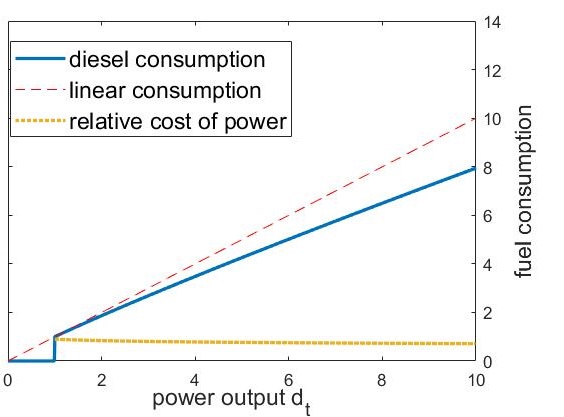}
            \caption{$\rho(d)=d^{0.9}$}    
        \end{subfigure}
\caption{The panels above show two examples of efficiency function (litres/KW), on the left $\rho(d)=\frac{(d-6)^3+6^3+d}{10}$, typical of a generator designed to operate at medium regime, on the right $\rho(d)=d^{0.9}$, typical of a generator designed to operate a full capacity.}
\label{fig:efficiency_diesel}   
    \end{figure}

Notice that, in addition, when the engine is turned ON, an extra amount of fuel is burned in order for the generator to warm up and reach working regime. We model the cost of burning extra fuel with a switching cost $\CK$ that is paid every time the switch changes from $0$ to $1$. The fuel consumption of the diesel generator is modeled by an increasing function $\rho(d_t)$ which maps the power $d_t$ produced during one time step into the quantity of diesel necessary for such output. Denoting by $P_t$ the price of fuel at time $t$, the cost of producing $d_t$ KW of power at one time step is $P_t\rho(d_t)$; for simplicity we take a constant price of the fuel $P_t=p$.  Two examples of efficiency functions $\rho$ are described in figure \ref{fig:efficiency_diesel}.

%\begin{figure}[h]
%\includegraphics[width=0.6\textwidth]{efficiency}
%\caption{The panels above show two examples of efficiency function, on the left we find $\rho(d)=\frac{(d-6)^3+6^3+d}{10}$ which represents a generator which best efficiency is producing at medium regime, on the right $\rho(d)=d^{0.9}$ which represents a generator whose best efficiency is attained at maximum regime.}
%\label{fig:efficiency_diesel}
%\end{figure}

\subsection{Dynamics of the Battery}
  
The storage device is directly connected to the microgrid and therefore its output is equal to the imbalance between demand $X_t$ and diesel generator output $d_t$, when this is allowed by the physical constraint. The battery therefore is discharged in case of insufficiency of the diesel output and charged when the diesel generator and renewables provide a surplus of power.

Let us denote the power output of the battery by $B^d_t$ and its power rating by $B^{\max}$ and $B^{\min}$, where  $B^{\max}$ and $B^{\min}$ represent respectively  the maximum and minimum output. Thus:
\begin{equation}
\label{eq:battery_output}
B^d_{t} = \frac{I_{t}^d-I_{\max}}{\Delta t} \vee\big(B^{\min} \vee (X_{t} - d_{t}) \wedge B^{\max}\big)\wedge\frac{I_t^d}{\Delta t} 
\end{equation}
The case where $B^d_t <0$, represents that the battery is  charging while the case where $B^d_t>0$, represents that the battery is supplying power.

 Notice then that an energy storage has a limited amount of capacity after which it can not be charged further, as well as an ``empty'' level below which no more power can be provided from the battery. We denote the state of charge by the controlled process $I^d_t$ which is described by the following equation:
\begin{equation}
\label{eq:inventory_dynamics}
I^d_{t+1}=I^d_t-B^d_{t}\Delta t,\quad t\in\{0,1,\dots ,T-1\},\quad I^d_0=w_0
\end{equation}
here $I^d_t\in[0,I_{max}]$ and $B^d_t\in[B^{\min},B^{\max}]$, for $B^{\min}<0$ and $B^{\max}>0$. 
For simplicity we assume that the battery is $100\%$ efficient. Notice that we used superscript $d$ on $B^d$ and $I^d$ to highlight the dependence of these processes on the controlled diesel output $d_t$.

Intuition tells us that the bigger the battery, the less diesel will be needed to run the operations of the microgrid. This is true because a bigger battery would allow to store for later use a bigger proportion of the excess power produced by the renewables. Batteries however are very expensive, and the cost per KWh of capacity scales almost linearly for the kind of devices we consider in this paper (parallel connection of smaller batteries), hence it is important  to find the optimal size of battery for the needs of each specific microgrid. 

\subsection{Management of the Microgrid}
\label{microgrid_management}

The purpose of the microgrid is to provide a cheap and reliable source of power supply to at least match the demand. Therefore, we search for a control policy for the diesel generator which minimizes the operating cost and produces enough electricity to match the residual demand.
In order to assess how well we are doing in supplying electricity, we introduce the controlled imbalance process $S_t$ defined as follows: 
\begin{equation}\label{st}
S_{t}=X_{t}-B^d_{t}-d_t \quad t\in[0,T]
\end{equation}
Ideally, the owner of the Microgrid would like to have $S_t =0 \quad \forall \; t$. This situation represents the perfect balance of demand and generation. When $S_{t}>0$ we observe a \textit{blackout}, residual demand is greater than the production meaning that some loads are automatically disconnected from the system. The situation $S_t<0$ is defined as a curtailment of renewable resources and takes place when we have a surplus of electricity. 

We treat the two scenarios, blackout and curtailment asymmetrically.  To ensure no-blackout $S_t \leq 0$ and regular supply of power, we impose a constraint on the set of admissible controls:
\begin{equation}
\label{eq:blackout_constraint}
\begin{split}
&S_t\le0 \\
\text{i.e. } &d_t \ge X_{t}-B^d_{t}.
\end{split}
\end{equation}

However, for $S_t<0$ i.e. surplus of electricity, we penalize the microgrid using a proportional cost denoted by $C$. Large penalty would lead to low level of curtailment and can be thought of as a parameter in the subsequent optimization problem.

% The blackout is an undesirable situation and to avoid it we impose a constraint, on the set of admissible controls:

% Given however the reduced size of the microgrid and the intent to maintain the size of the battery small, we should allow the imbalance process to become negative; the situation $S^d_t <0$ is defined as a curtailment of renewable resources and takes place when we have a surplus of electricity. In order to have control over the quantity of curtailed energy resulting from a given control policy for the diesel generator, we introduce a proportional cost denoted by $C$. This parameter can be tuned in the optimization problem to achieve the desired level of curtailment

A rigorous mathematical description of the microgrid management problem follows in section \ref{control_problem}.

%%%%%%%%%%%%%%%%%%%%%%%%%%%%%%%%%%%%%%%%%%%%%%%%%%%%%%%%%%%%%%%%%%%%%%%%%%%%%%%%%%%%%%%%%%%%%%%%%%%%%%%%%%%%%%%%%%%%%%%%%%%%%%%%%%%%%%%%%%%%%%%%%%%%%%%%%%%%%%%%%%%%%%%%%%%%%%%%%%%%%%%%%%%%%%%%%%%%%%%%%

\section{Stochastic optimization problem}
\label{control_problem}
We state now the stochastic control problem for the diesel generator operating in a microgrid system as described in section \ref{model_description}. In practice we seek a control that minimizes the cost of diesel usage $p \rho(d)$, the switching cost $\CK$ and the curtailment cost $C |S_t|\mathds{1}_{\{S_t<0\}}$, under the no black-out constraint $S_t\le0$. 

Note that, given the type of control we have on the diesel generator, we can frame the optimization problem as a special case of stochastic control problems known as  optimal switching problems.

Let us denote  by $\mathcal{F}_t$ the filtration generated by the residual demand process $(X_s)_{s=0}^t$, the state of charge process $(I_s^d)_{s=0}^t$ and the current regime $m_t$, which represents all the information available on the system up to time $t$. In practice, given the markovianity of the problem, we have that $\mathcal{F}_t$ is reduced to the $\sigma$-field generated by the triple $(X_t,I_t^d,m_t)$.

Let us define the pathwise value $\mathcal{J}$, given by
%\begin{equation}
%\label{eq:pathwise_value}
%\mathcal{J}^{m}(t,\textbf{X}_{t},\textbf{I}_{t}; d_{t},m_{t})=\sum_{s=t}^{T-1}\mathds{1}_{\{m_{s+1}-m_{s}=1 \}}K+p\rho(d_s)+C S_s\mathds{1}_{\{S_s<0\}}+g(I^{d}_T).
%\end{equation}
%where $(\textbf{X}_{t},\textbf{I}_{t}^d; d_{t},m_{t}) = ( X_s, I^{d}_s; d_{s},m_{s} )_{s=t}^{T}$.
%
\begin{equation}
\label{eq:pathwise_value}
\mathcal{J}(t,\textbf{X}_{t},\textbf{I}_{t},m_t; d_{t})=\sum_{s=t}^{T-1}\mathds{1}_{\{m_{s+1}-m_{s}=1 \}}\CK+p\rho(d_s)+C|S_s|\mathds{1}_{\{S_s<0\}}+g(I^{d}_T).
\end{equation}
%where $(\textbf{X}_{t},\textbf{I}_{t}) = ( X_s, I^{d}_s )_{s=t}^{T}$.
where $(\textbf{X}_{t},\textbf{I}_{t}, m_t; d_{t}) = ( X_s, I^{d}_s, m_s; d_{s} )_{s=t}^{T}$.
 As a consequence, we define the value function as:
\begin{equation}
\label{eq:value_function}
V(t,x,w, m)=\min_{d_{t}=(d_{u})_{u=t}^T}\left\{\mathbb{E}\left[ \mathcal{J}(t,\textbf{X}_{t},\textbf{I}^{d}_{t},m_{t}; d_{t})\Big| X_t=x, I^{d}_t=w ,m_t=m \right] \right\}
\end{equation}
\begin{subequations}
\begin{align}
\text{subject to }\quad & d_t\ge X_{t}-B^d_{t} \quad\;\forall t \label{cond_st} \\
&  d_t \in [d_{min},d_{max}]\cup \{0\}. \label{cond_dt}\\
& B^d_{t} = \frac{I_{t}^d-I_{\max}}{\Delta t} \vee\big(B^{\min} \vee (X_{t} - d_{t}) \wedge B^{\max}\big)\wedge\frac{I_t^d}{\Delta t} \label{cond_bat}
\end{align}
 \end{subequations}
where \eqref{cond_st} represents the black-out constraints translated for the power produced by the diesel generator, \eqref{cond_dt} represents the minimum and maximum power output of the generator and \eqref{cond_bat} models the physical constraints of the battery: maximum input/output power and maximum capacity.

From equation \eqref{eq:value_function}, we can write the associated dynamic programming formulation which helps understand the structure of the problem composed of two optimal control problems: an optimal switching problem between being in the regime ON or OFF, and another absolutely continuous control problem assuming the regime is ON. The equation reads as follows:
% \begin{equation}
% \label{eq:dpe}
% V(t,x,w,m)= \min\Big\{\mathcal{V}(t,x,w,m;0,0),\underset{d\in\mathcal{U}_{t}}{\min}\Big(\mathcal{V}(t,x,w,m;d,1)\Big)+K\mathds{1}_{\{m=0\}}\Big\},
% \end{equation}
% where 
% \[
% \mathcal{V}(t,x,w,m;d,y)=p\rho(d_{t})+C|S_t|\mathds{1}_{\{S_t<0\}} + \bE\big[ V(t+1,X_{t+1},I^{d}_{t+1},m_{t+1}=y)\big|X_t=x,I^d_t=w,m_t=m \big],
% \]
% and $\mathcal{U}_{t}$ is the collection of admissible controls $d$ at each time step $t$,i.e. 
% \begin{equation}
% \mathcal{U}_{t}:=\{d_t: \text{equations \eqref{cond_st} - \eqref{cond_bat} are satisfied and $d_t$ adapted to } \mathcal{F}_t\}.
% \end{equation}

%%%%%%%%%%%%%%
% behaviour of m_t is similar to I_t, both of them are F_{t-1} measurable. As a result, $m_{t+1}$ is decided because of the control chosen at time $t$. Yes, I was using a fancy definition that apparently was misleading.
%%%%%%%%%%

%\begin{equation}
%\label{eq:dpe_old}
%V(t,x,w,m)= \min\Big\{\mathcal{V}(t,x,w;0),\underset{d\in\mathcal{U}_{t}\setminus \{0 \}}{\min}\Big(\mathcal{V}(t,x,w;d)\Big)+\CK\mathds{1}_{\{m=0\}}\Big\},
%\end{equation}
%where 
%\[
%\mathcal{V}(t,x,w;d)=p\rho(d_{t})+C|S_t|\mathds{1}_{\{S_t<0\}} + \bE\big[ V(t+1,X_{t+1},I^{d}_{t+1},m_{t+1}=\mathds{1}_{d>0})\big|X_t=x,I^d_t=w,d_t=d \big],
%\]

\begin{equation}
\label{eq:dpe}
V(t,x,w,m)= \underset{d\in\mathcal{U}_{t}}{\min}\Big(\mathds{1}_{\{m_{t+1}-m_{t}=1 \}}\CK + p\rho(d)+C|S_t|\mathds{1}_{\{S_t<0\}} + \mathcal{C}(t,x,w,m;d)\Big),
\end{equation}
where 
\[
\mathcal{C}(t,x,w,m;d) =  \mathbb{E}[V(t+1,X_{t+1},I_{t+1},m_{t+1})|X_t=x,I_t=w,d_t=d,m_t=m],
\]
is the conditional expectation of the future costs and $\mathcal{U}_{t}$ is the collection of admissible controls $d$ at each time step $t$, i.e. 
\begin{equation}
\mathcal{U}_{t}:=\{d_t: \text{equations \eqref{cond_st} - \eqref{cond_bat} are satisfied and $d_t$ adapted to } \mathcal{F}_t\}.
\end{equation}

In order to ensure that the set of admissible controls is nonempty we introduce the following assumption:

\begin{assumption}
\label{assumption:diesel_power}
The diesel generator is powerful enough to supply demand at all times, i.e there is always a control $d$ that satisfies the blackout constraint.
\end{assumption}
\begin{remark}
\label{r:truncation}
We enforce assumption \ref{assumption:diesel_power} by redefining the residual demand process with a truncated version of \eqref{eq:residual_demand}, such that $\tilde{X}_t = \min (X_t, X_{\max}) $ is the residual demand. In practice this is reasonable because the  maximum power that could be required from the microgrid is known apriori and the diesel generator is generally sized to the maximum capacity installed on the system.  For the sake of notational simplicity, we will drop the  $\sim$ on the variable $\tilde{X}_t$  from the following sections.
\end{remark}

Note that \eqref{eq:dpe} provides a direct technique to solve problem \eqref{eq:value_function}, iterating backward in time from a known terminal condition and solving a static, one period, optimization problem at each time step. The only difficulty in this procedure lies in the estimation of conditional expectations of future value function, which can not be computed exactly. In the next section \ref{numerical_methods} we will focus on the numerical solution of \eqref{eq:value_function}.

%%%%%%%%%%%%%%%%%%%%%%%%%%%%%%%%%%%%%%%%%%%%%%%%%%%%%%%%%%%%%%%%%%%%%%%%%%%%%%%%%%%%%%%%%%%%%%%%%%%%%%%%%%%%%%%%%%%%%%%%%%%%%%%%%%%%%%%%%%%%%%%%%%%%%%%%%%%%%%%%%%%%%%%%%%%%%%%%%%%%%%%%%%%%%%%%%%%%%%%%%

\section{Numerical Resolution}
\label{numerical_methods}

In this section we describe the algorithm which we want to employ in the solution of the energy management problem for the  Microgrid system described in section \ref{control_problem}. The main mathematical difficulty comes from the approximation of conditional expectations in \eqref{eq:dpe}, which we will tackle using a family of methods called Regression Monte Carlo.

The algorithm we propose fully exploits the dynamic programming formulation \eqref{eq:dpe}: we start generating a set of simulations (scenarios) of the process $X$, which we will refer to as training points, then we optimize our policy so that it performs well, on average (weighted on the probability of each scenario), on the different scenarios. 

In practice, we initialize the value function at last time step in the backward procedure to be equal to the terminal condition $g$. We then iterate backward in time and at each time step over each training point we choose the control that minimizes the sum of one step cost function and the estimated conditional expectation of the future costs $\tilde{\mathcal{C}}(t,x,w,m;d)$. Note that, as expected, the conditional expectation is a function of time, the state of the system $(x,w)$ and the state of the diesel generator, represented by the ON/OFF switch $m$ and the control $d$.

As the iteration reaches the initial time point we collect a set of optimal actions for each time step and many different scenarios; in addition, since the problem is Markovian, we can summarize such strategies in the form of control maps: best action at each time $t$ given a pair of state variables $(X_t,I_t)$ and state of the diesel generator $m_t$. We propose three different techniques to compute $\tilde{\mathcal{C}}$ in section \ref{regression_c}. 
 
A fair assessment of the quality of the control policies approximated by the algorithm just introduced is obtained by running a number of forward Monte Carlo simulations of the residual demand, controlling the system using such policies and then taking the average performance.

We give a general description of the pseudo code in algorithm \ref{algo:microgrid_management}. 

\begin{remark}
Notice that it is typical of Regression Monte Carlo algorithms to provide the optimal policy only implicitly, in the form of minimizer of an explicit parameterized function. The outputs of the algorithm are therefore the parameters (regression coefficients) of such function. 
\end{remark}

\begin{algorithm}
\caption{Regression Monte Carlo algorithm for Microgrid management}
\textbf{input:} number of basis $K$, number of training points $M$, discretisation of the inventory $D$, time-steps $N$.
\begin{algorithmic}[1]
\State \textbf{optimization:}
\If{Inventory discretisation}
\State Generate a customary grid $\{w_0,\ldots,w_D \}$ points over the domain of $I_t$.
\State Simulate $\{X_{t}^{j}\}_{j,t=1}^{M',N}$ according to its dynamics where $M' = M/(D+1)$;
\State Define $\{X_{t}^{j},I_{t}^{j}\}_{j=1}^{M}$ as cross product of  $\{X_{t}^{j}\}_{j=1}^{M'}$ and $\{w_j\}_{j=0}^{D}$ for $\forall t$
\EndIf  
\If{Regression 2D}
\If{Regress Later}
\State  Generate $\{X_{t}^{j},I_{t}^{j}\}_{j,t=1}^{M,N}$ accordingly to a distribution $\mu$;
\EndIf
\If{Regress Now}
\State  Generate $\{X_{t}^{j}\}_{j,t=1}^{M,N}$ according to its dynamics and  $\{I_{t}^{j}\}_{j,t=1}^{M,N}$ according to a distribution $\mu$;
\EndIf
\EndIf  
\State Initialize the value function $V(N, X_{N}^{j},I_{N}^j,1) =V(N, X_{N}^{j},I_{N}^j,0) =g(I_N^j), \quad \forall j=1,\,\dots,\, M$;
\For{$t=N$ to $1$}
\State Compute the approximated continuation value $\tilde{\mathcal{C}}$ using Algorithms \ref{algo:regression_2D} or \ref{algo:regression_1D}
\For{$j=1$ to $M$}
\For{$m=0$ to $1$}
\State $F = \tilde{\mathcal{C}}(X_{t}^j,I_{t}^j; 0, 0)$
\State 
\[
V(t, X_{t}^{j},I_{t}^j,m)=
\begin{cases}
\Big( \min\limits_{d\in\mathcal{U}_{t} \setminus \{0 \} }\Big\{ p\rho(d)+ C|S_t|\mathds{1}_{\{S_t<0\}} + \tilde{\mathcal{C}}(X_{t}^j,I_{t}^j;1,d) \Big\}+\CK\mathds{1}_{\{m=0\}} \Big) \wedge F  &\quad \text{if } 0 \in \mathcal{U}_{t}\\ 
\min\limits_{d\in\mathcal{U}_{t} }\Big\{ p\rho(d)+ C|S_t|\mathds{1}_{\{S_t<0\}} + \tilde{\mathcal{C}}(X_{t}^j,I_{t}^j;1,d) \Big\}+\CK\mathds{1}_{\{m=0\}} &\quad \text{otherwise}
\end{cases}
\]
\EndFor
\EndFor
\EndFor
\State \textbf{simulation:}
\State initialize processes
\For{$t=1$ to $N-1$}
\For{$j=1$ to $M$}
\State $F_1 = \tilde{\mathcal{C}}(X_{t}^j,I_{t}^j;0, 0)$
\State $F_2 = \min\limits_{d\in\mathcal{U}_{t} \setminus \{0 \} }\Big\{ p\rho(d) +C|S_t|\mathds{1}_{\{S_t<0\}} +\tilde{\mathcal{C}}(X_{t}^j,I_{t}^j;1,d) \Big\}+\CK\mathds{1}_{\{m_{t}^j=0\}}$
\State $m_{t+1}^j= \mathds{1}_{\{ (0 \notin \mathcal{U}_{t} ) \text{ or }  (0 \in \mathcal{U}_{t} \text{ and } F_2<F_1)\} }$

\If{$m_{t+1}^j=1$}
\State $d_{t}=\argmin\limits_{d\in\mathcal{U}_{t}}\Big\{ p\rho(d) +C|S_t|\mathds{1}_{\{S_t<0\}} +\tilde{\mathcal{C}}(X_{t}^j,I_{t}^j;1,d) \Big\}$
\EndIf
\State compute $X_{t+1}^j$ and $I_{t+1}^j=I_{t}^j-B^d_{t}\Delta t$
\State $J_{t+1}^j=J_{t}^j+p\rho(d_{t}) + C|S_t|\mathds{1}_{\{S_t<0\}}+\CK\mathds{1}_{\{m_{t+1}- m_t=1\}}$
\EndFor
\EndFor
\State $V(0,x,w, m)=\frac{1}{M}\sum_{j=1}^M (J_{N}^j+g(I_N^j))$
\end{algorithmic}
\textbf{output:} control policy $\{d_t\}$, value function $V$.
\label{algo:microgrid_management}
\end{algorithm}

\subsection{Regression for continuation value}%%%%%%%%%%%%%%%%%%%%%%%%%%%%%%%%%%%%%%%%%%%%%%
\label{regression_c}

In this section we present the numerical techniques we use to estimate conditional expectations $\mathcal{C}(t,x,w,m;d)$ in algorithm \ref{algo:microgrid_management}. These techniques belong to the realm of Regression Monte Carlo methods, and in particular these specifications allow to deal with degenerate controlled processes (the inventory). We focus on two main variants: a two dimensional approximation of the conditional expectation and a discretisation technique which considers a collection of one dimensional approximations.

% In particular, we test three algorithms: Grid Discretization, which is characterized by a one dimensional "Regress Now" projection in the residual demand dimension repeated at different inventory points and 2D Regression, which uses a two dimensional regression in residual demand and inventory. The 2D Regression uses a "Regress Later" technique in the inventory dimension, while we test both "Regress Now" and "Regress Later" in the other dimension. With the terms ``Now'' and ``Later'' we indicate projections over basis functions measurable with respect to $\mathcal{F}_t$ and $\mathcal{F}_{t+1}$ respectively. For details on these techniques see \citet{Balata2017i},\citet{?} or \citet{?}. Note that in the three algorithm we repeat the regression approximation for both values of $m$. 
 
%Let us denote by $\{X_t^m\}_{m=1}^M$ the collection of training points at time $t$, similar notation is used for the inventory $\{I_t^m\}_{m=1}^M$.

In particular, we test three algorithms: Grid Discretisation, Regress Now and Regress Later. Grid Discretization is characterized by a one dimensional projection in the residual demand dimension repeated at different inventory points. Regress Now/Later, on the other hand, use a two dimensional regression in residual demand and inventory. Moreover, while Grid Discretization and Regress Now require projection of the value function at $t+1$ on $\mathcal{F}_t$ measurable basis functions, Regress Later requires an $\mathcal{F}_{t+1}$ projection. For details on these techniques see \citet{balata17} for regress later,\citet{boogert08,warin} for GD and \citet{ludkovski10} for 2D regress now.  Note that in the three algorithms we repeat the regression approximation for both values of $m$. An open source platform has also been developed  to numerically solve wide variety of stochastic optimization problems in \citet{StOpt}. 
 
Let us denote by $\{X_t^j\}_{j=1}^M$ the collection of training points at time $t$, similar notation is used for the inventory $\{I_t^j\}_{j=1}^M$.

\subsubsection{Grid Discretisation}
Grid discretisation is characterized by a one dimensional approximation of the conditional expectation repeated at different levels of inventory.
Let $\Upsilon_I=\{w_0=0,\dots,w_D=I_{max}\}$ be a discretisation of the state space of the inventory and $\{X_t^j\}_{j=1,t=1}^{M,N}$ be generated from a forward simulation of the dynamics of $X$. We define the approximation of the continuation value on the grid $\Upsilon_I$ by regressing the set of value functions  $\{V(t+1,X_{t+1}^j,w_i) \}_{j=1}^{M}$ over the basis functions $\{\phi_k(x)\}_{k=1}^K$ for each $\{w_i\}_{i=0}^D$,  obtaining:
\[
\hat{\mathcal{C}}(t,x,w_i;m)=\sum_{k=1}^K \alpha_{k,i,m}^t\phi_k(x)\,,\quad i=0,\,1,\,\dots,\,D,
\]
where we compute a collection of regression coefficients through least square minimization
\[
\boldsymbol{\alpha}^{t}_{i,m}=\argmin\limits_{a\in\mathbb{R}^K}\Big\{\frac{1}{M}\sum_{j=1}^M\big(V(t+1,X_{t+1}^j,w_i,m)-\sum_{k=1}^K a_{k}\phi(X_{t}^j)\big)^2\Big\},
\]
where we define $\mathbb{R}^K\ni\boldsymbol{\alpha}^{t}_{i,m}=(\alpha^{t}_{1,i,m},\,\dots,\,\alpha^{t}_{K,i,m})$.

Note that the least square projection is a sample estimation of the $L^2$ projection induced by the conditional expectation, for this reason we can approximate the function $\mathcal{C}(t,\cdot)$ using a least square projection of the value function at time $t+1$.
%However, as we have not included the inventory in the basis functions, we need to interpolate between values of $\hat{\mathcal{C}}(t,x,w_j;m)$ in order to obtain an estimation of the value function for $I_t\in(w_j,w_{j+1})$. Let us define by $\mathcal{C}(t,x,w;m,d)$ the linear interpolation
%where $\omega(t,w,d)=\frac{w_{j+1}-w-B_t^d}{w_{j+1}-w_j}$ and $j=0,\,\dots,\,D-1$.
However, as we have not included the inventory in the basis functions, we need to interpolate between values of $\hat{\mathcal{C}}(t,x,w_i;m)$ in order to obtain an estimation of the value function for $I_t\in(w_i,w_{i+1})$. Let us define by $\tilde{\mathcal{C}}(t,x,w;m,d)$ the linear interpolation
\[
\tilde{\mathcal{C}}(t,x,w;m,d)=\omega(t,w,d)\hat{\mathcal{C}}(t,x,w_i,m)+\big(1-\omega(t,w,d)\big)\hat{\mathcal{C}}(t,x,w_{i+1},m)\,,\quad w-B_t^d\Delta t\in[w_i,w_{i+1}),
\]
where $\omega(t,w,d)=\frac{w_{i+1}-w+B_t^d\Delta t}{w_{i+1}-w_i}$ and $i=0,\,\dots,\,D$.\\
Details of the algorithms are given in the pseudocode  \ref{algo:regression_1D}.
 \begin{algorithm}
 \caption{Regression technique for continuation value: Grid Discretisation}
 \textbf{input:} $\{V(t+1,X_{t+1}^j,I_{t+1}^j,m) \}_{j=1}^M$, $\{\phi_k\}_{k=1}^K$.
 \begin{algorithmic}[1]
 \For{$i=0$ to $D$}
  \State  $  \boldsymbol{\alpha}^{t}_{m}=\argmin\limits_{a}\Big\{\sum\limits_{j=1}^M \Big(V(t+1, X_{t+1}^{j},w_i, m)-\sum\limits_{k=1}^K a_k\phi_k(X^j_{t})  \Big)^2\Big\}$;
    \State Define $\hat{\mathcal{C}}(t,x,w_i,m)=\sum_{k=1}^K\alpha^{t}_{k,i,m} \phi_k(x)$, $m=0,1$;
    \EndFor
    %  \State Define $\mathcal{C}(t,x,w;m,d)=\frac{w_{j+1}-w-B_t^d}{w_{j+1}-w_j}\hat{\mathcal{C}}(t,x,w_j;m,d)+\frac{w+B_t^d-w_j}{w_{j+1}-w_j}\hat{\mathcal{C}}(t,x,w_{j+1};m,d)$, $w\in[w_j,w_{j+1})$, $m=0,1$.
  \State Define $\tilde{\mathcal{C}}(t,x,w;m,d)=\frac{w_{i+1}-w+B_t^d\Delta t}{w_{i+1}-w_i}\hat{\mathcal{C}}(t,x,w_i;m,d)+\frac{w-B_t^d-w_i}{w_{i+1}-w_i}\hat{\mathcal{C}}(t,x,w_{i+1};m,d)$, $w\in[w_i,w_{i+1})$, $m=0,1$.
\end{algorithmic}
\textbf{output:} $\tilde{\mathcal{C}}$, $\{\alpha^{t}_{k,i,m}\}_{k=1,i=1,m=0}^{K,D,1}$.
 \label{algo:regression_1D}
 \end{algorithm}
\subsubsection{2D Regression}

Contrary to the grid discretisation approach, the 2D regression methods approximate the conditional expectation of the value function as a surface, function of both residual demand $X$ and inventory $I$, without the need for interpolation. 
%Such approximation is made possible by the so-called ``Regress Later'' approximation, which allow to extend the Regression Monte Carlo approach to controlled processes by projecting the value function at time $t+1$ on basis functions $\{\phi_k(\cdot)\}_{k=1}^K$ that are $\mathcal{F}_{t+1}$ measurable.
 In the problem we consider, the control only acts on a degenerate (deterministic) process and we can therefore test two specifications of the method: ``Regress Now'', where we project over  $\{\phi_k(X_t,I_{t+1})\}_{k=1}^K$ and ``Regress Later'', where we project over  $\{\phi_k(X_{t+1},I_{t+1})\}_{k=1}^K$. The terminology Regress Now or Regress Later is attributed to the time step of the exogenous variable $X_t$ used in the projection.
 
 In Regress Now, we generate training points $\{X_t^j\}_{j=1,t=1}^{M,N}$ from a forward simulation of the dynamics of $X$ and $\{I_t^j\}_{j=1,t=1}^{M,N}$ from a distribution $\mu_N$ on $[0,I_{max}]$. In Regress Later, on the other hand, we generate both processes $\{X_t^j,I_t^j\}_{j=1,t=1}^{M,N}$ from an appropriate distribution $\mu_L$, for details see \citet{balata17}. In the following we will generalize the discussion of the two approaches by using the subscript $r$ with realization $t$ to indicate Regress Now algorithm and $t+1$ to indicate Regress Later.  As training measures we choose $\mu_N$ to be the Lebesgue measure on $[0,I_{\max}]$ and $\mu_L$ to be Lesbegue measure on $[0,I_{\max}]\times[-X_{\max}, X_{\max}]$.
 
 The regression coefficients in the 2D regression Monte Carlo method are computed by least-square projection as:
 \[
 \boldsymbol{\alpha}^{t}_{m}=\argmin\limits_{a\in\mathbb{R}^K}\Big\{\frac{1}{M}\sum_{j=1}^M\big(V(t+1,X_{t+1}^j,I_{t+1}^j,m)-\sum_{k=1}^K a_k\phi(X_{r}^j,I_{t+1}^j)\big)^2\Big\},
\]
where we define $\mathbb{R}^K\ni\boldsymbol{\alpha}^{t}_{m}=(\alpha^{t}_{1,m},\,\dots,\,\alpha^{t}_{K,m})$.

Let us recall, denoting by $\boldsymbol{\phi}$ the vector $\big(\phi_1(\cdot),\,\dots,\,\phi_K(\cdot)\big)$, that the coefficients $ \boldsymbol{\alpha}^{t}_{m}$  can be computed explicitly by 
\[
 \boldsymbol{\alpha}^{t}_{m}=\Big(\mathbb{E}_{\mu}\big[\boldsymbol{\phi} \boldsymbol{\phi}^T\big] \Big)^{-1} \mathbb{E}_{\mu}\Big[V(t+1,X_{t+1},I_{t+1},m)\boldsymbol{\phi}\Big]^T\approx \Big(\sum_{j=1}^M \boldsymbol{\phi} \boldsymbol{\phi}^T \Big)^{-1} \sum_{j=1}^M V(t+1,X_{t+1}^j,I_{t+1}^j,m)\boldsymbol{\phi}^T
\]
and therefore, even though the regression coefficients are random (sample average approximation of expectations with respect to the measure $\mu$) they are independent of $\mathcal{F}_{t}$. Given the previous remark we can estimate the conditional expectation of future value through:
\[
\tilde{\mathcal{C}}(t,x,w;m,d)=\mathbb{E}\Big[\sum_{k=1}^K\alpha_{k,m}^t\phi_k(X_{r},I_{t+1})\Big|\mathcal{F}_{t}\Big]=\sum_{k=1}^K\alpha_{k,m}^t\mathbb{E}\Big[\phi_k(X_{r},I_{t+1})\Big|X_t=x,\,I_t=w,\,d_t=d \Big].
\]

The explicit value of $\mathbb{E}\Big[\phi_k(X_{r},I_{t+1})\Big|X_t=x,\,I_t=w,\,d_t=d\Big]$ now depends on $r$, i.e. whether we are using ``Regress Now'' or ``Regress Later'' to deal with the uncontrolled residual demand. In the first case we simply obtain, from the measurability of $X_t$, 
%\[
%\mathbb{E}\Big[\phi_k(X_{t},I_{t+1})\Big|\mathcal{F}_{t}\Big]=\phi_k(x,w+B^d_t)=:\tilde{\phi}_k(x,w,d).
%\]
\[
\mathbb{E}\Big[\phi_k(X_{t},I_{t+1})\Big|\mathcal{F}_{t}\Big]=\phi_k(x,w-B^d_t\Delta t)=:\tilde{\phi}_k(x,w,d).
\]
 In the second case we need to compute the  expectation with respect to the randomness contained in the transition function from $X_t$ to $X_{t+1}$ and we simply write 
% \[
% \mathbb{E}\Big[\phi_k(X_{t+1},I_{t+1})\Big|\mathcal{F}_{t}\Big]=\mathbb{E}_{\xi}\Big[\phi_k(x+ b(\Lambda_{t} -x )\Delta t + \sigma \sqrt{\Delta t} \xi ,w+B^d_t)\Big]=:\hat{\phi}_k(x,w,d).
% \]
\[
 \mathbb{E}\Big[\phi_k(X_{t+1},I_{t+1})\Big|\mathcal{F}_{t}\Big]=\mathbb{E}_{\xi}\Big[\phi_k(x+ b(\Lambda_{t} -x )\Delta t + \sigma \sqrt{\Delta t} \xi ,w-B^d_t\Delta t)\Big]=:\hat{\phi}_k(x,w,d).
 \]
 
 \begin{remark}
 For polynomial basis functions, i.e. $\phi_k(X_{t+1},I_{t+1}):=X^p_{t+1}I^q_{t+1}$, the conditional expectation $\hat{\phi}_k(x,w,d)$ can be written in closed form as:
 \[
 \begin{split}
 \hat{\phi}_k(x,w,d)
& =\mathbb{E}\big[X^p_{t+1},I^q_{t+1}\big|X_t=x,I_t=w,d_t=d\big]\\
& =I_{t+1}^q\sigma^p dt^{\frac{p}{2}}\sum_{k=0}^p\mathbb{I}_{\{(p-k)\text{ is odd}\}}{p \choose k}\Big(x\frac{1-\lambda dt}{\sigma\sqrt{dt}}\Big)^k\prod_{j=1}^{\frac{p-k}{2}}(2j-1)
\end{split}
 \] 
 \end{remark}
  Using the notation just introduced we can summarize the differences between the two techniques in the following table:

%\begin{table}[H]
%\centering
%\begin{tabular}{l|l|l|l}
%   & $\phi_k$              & $\mathbb{E}[\phi_k|X_t,I_t,d_t]$                    & $\mathcal{C}(t,x,w;d,m)$             \\ \hline\hline
%RN & $(X_t,I_{t+1})$     & $\phi_k(X_t,I_t+B^d_t)$                 & $\sum_{k=1}^K\alpha_k^t\hat{\phi}_k(x,w,d)$ \\ \vspace{5pt}
%RL & $(X_{t+1},I_{t+1})$ & $\mathbb{E}[\phi_k(X_t+ b(\Lambda_{t} -X_t )\Delta t + \sigma \sqrt{\Delta t} \xi ,I_t+B^d_t)]$ & $\sum_{k=1}^K\alpha_k^t\tilde{\phi}_k(x,w,d)$                                
%\end{tabular}
%\end{table}
\begin{table}[H]
\centering
\begin{tabular}{l|l|l|l}
   & $\phi_k$              & $\mathbb{E}[\phi_k|X_t,I_t,d_t]$                    & $\mathcal{C}(t,x,w,m;d)$             \\ \hline\hline
RN & $(X_t,I_{t+1})$     & $\phi_k(X_t,I_t-B^d_t\Delta t)$                 & $\sum_{k=1}^K\alpha_{k,m}^t\tilde{\phi}_k(x,w,d)$ \\ \vspace{5pt}
RL & $(X_{t+1},I_{t+1})$ & $\mathbb{E}[\phi_k(X_t+ b(\Lambda_{t} -X_t )\Delta t + \sigma \sqrt{\Delta t} \xi ,I_t-B^d_t\Delta t)]$ & $\sum_{k=1}^K\alpha_{k,m}^t\hat{\phi}_k(x,w,d)$                               
\end{tabular}
\end{table}
Details of the algorithms are given in the pseudocode \ref{algo:regression_2D} .

\begin{algorithm}
 \caption{Regression technique for continuation value: 2D Regression}
 \textbf{input:} $\{V(t+1,X_{t+1}^j,I_{t+1}^j,m) \}_{j=1}^M$, $\{\phi_k\}_{k=1}^K$.
 \begin{algorithmic}[1]
 \If {Regress Later}
 \State $r=t+1$
 \ElsIf {Regress Now}
   \State $r=t$
   \EndIf
  \State  $ \boldsymbol{\alpha}^{t}_{m}=\argmin\limits_{a}\Big\{\sum\limits_{j=1}^M \Big(V(t+1, X_{r}^{j},I_{t+1}^j,m)-\sum\limits_{k=1}^K a_k\phi_k(X^j_{r},I^j_{t+1})  \Big)^2\Big\}$, $m=0,1$;
  \State Define $\tilde{\mathcal{C}}(t,x,w;m,d)=\sum_{k=1}^K\alpha^{t}_{k,m} \mathbb{E}[\phi_k(X_{r},I_{t+1})|x,w,d]$ 
\end{algorithmic}
\textbf{output:} $\tilde{\mathcal{C}}$, $\{\alpha^{t}_{k,m}\}_{k=1,m=0}^{K,1}$.
\label{algo:regression_2D}
 \end{algorithm}

%%%%%%%%%%%%%%%%%%%%%%%%%%%%%%%%%%%%%%%%%%%%%%%%%%%%%%%%%%%%%%%%%%%%%%%%%%%%%%%%%%%%%%%%%%%%%%%%%%%%%%%%%%%%%%%%%%%%%%%%%%%%%%%%%%%%%%%%%%%%%%%%%%%%%%%%%%%%%%%%%%%%%%%%%%%%%%%%%%%%%%%%%%%%%%%%%%%%%%%%%

\section{Numerical Experiments}
\label{experiments}

In this section we use the algorithms introduced in section \ref{numerical_methods} to solve a simple instance of the microgrid management problem. We fix some base parameters and test the three algorithms; the one performing best is then used to study the sensitivity of the control policy and of the operational costs on changes in system parameters, hoping to gain some insight on the optimal design of the microgrid.

We now list the base parameters chosen for the numerical experiments; notice that the "s" column indicates whether a sensitivity analysis is run for such parameter. For the meaning of the parameters refer to section \ref{model_description}. 
\begin{table}[H]
\centering
    \begin{tabular}{|l||l|l|}
        \hline
        parameter & value & s \\ \hline
         $ T $ & $ 100 h $ &  \\      
         $ \Delta t $ & $ 0.25 h $ &  \\
	  $  b $ & $ 0.5 $ & *  \\
	  $ \sigma $ & $ 2 $ & * \\
	  $ \Lambda_t  $ & $ 0,\,\forall t $ &  \\
        \hline
    \end{tabular}\qquad     \begin{tabular}{|l||l|l|}
        \hline
        parameter & value & s \\ \hline
         $ I_{max} $ & $ 10 $ KWh & * \\
         $ \rho(d) $ & $ \frac{(d-d^*)^3+(d^*)^3+d}{10} \frac{litre}{\text{KW}}$ &  \\
	  $ d^*  $ & $ 6 $ KW& \\
	  $ p  $ & $ 1 $ \EUR & \\
	  $ g(i) $ & $ 0,\,\forall i $ &  \\
        \hline
    \end{tabular}\qquad     \begin{tabular}{|l||l|l|}
        \hline
        parameter & value & s \\ \hline
	  $ d_{min} $ & $ 1 $KW &  \\
	  $ d_{max} $ & $ 10 $KW &  \\
	  $ K $ & $ 5 $ \EUR& * \\
	  $ C $ & $ 0 $ \EUR& * \\
        \hline
    \end{tabular}
\end{table}

According to the parameters table above, and recalling remark \ref{r:truncation} the residual demand has the following dynamics:
\begin{equation}
\label{eq:base_parameters_X_dynamics}
X_{t+1}=  \big(  X_{t}(1-0.5 \Delta t) + \sigma \sqrt{\Delta t} \xi_{t} \big) \wedge 10, \quad t\in\{0,1,\dots,T-1\},
\end{equation}
where $\xi_t\sim\mathcal{N}(0,1)$.

We decided to use such simple dynamics for illustrative purposes in order to make the sensitivity of the optimal control policy to the remaining parameters more straight forward to understand.
 
\begin{figure}
        \centering
        \begin{subfigure}[b]{0.3\textwidth}
            \centering
            \includegraphics[width=\textwidth]{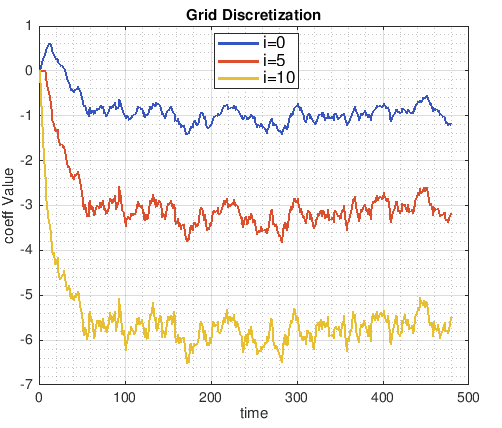}
            \caption[]%
            {{\small }}    
        \end{subfigure}
        \quad
        \begin{subfigure}[b]{0.3\textwidth}  
            \centering 
            \includegraphics[width=\textwidth]{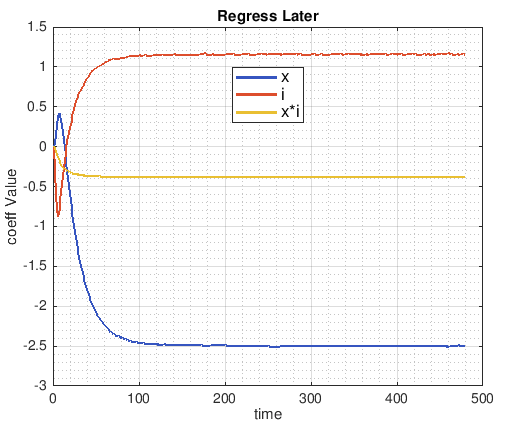}
            \caption[]%
            {{\small }}    
        \end{subfigure}
        \quad
        \begin{subfigure}[b]{0.3\textwidth}  
            \centering 
            \includegraphics[width=\textwidth]{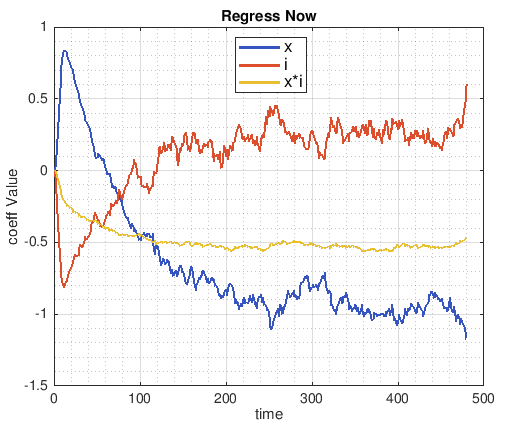}
            \caption[]%
            {{\small }}    
        \end{subfigure}
    \caption{In the three panels above we display the estimated regression coefficients corresponding to the basis $\{x,i,x\,i\}$ in the case of 2D regression, and $\{x\}$ at three different inventory levels for GD for $m_t=1$. Although we used basis function up to polynomial degree 2, we present few coefficients for clarity of presentation. Notice that the time axis is inverted to show the number of time steps computed backward. Remarkable smooth coefficients are computed by the Regress Later algorithm.}    
    \label{fig:stationary_coefficients}    
    \end{figure}

 Consider now that for the parameters listed above, the problem is time homogeneous. We have also observed empirically that the estimated continuation values tend to forget the terminal condition rather quickly. We show in Figure \ref{fig:stationary_coefficients} that the regression coefficients for all algorithms converge to a stationary value time steps, suggesting that optimization ran for longer time horizons would not bring any noticeable effect to control policy. Since all three methods use polynomial basis of degree two for the projection, it also allows for easy comparison of the dynamics of the coefficients across methods. For example, at inventory level $I=0$ the dynamics of the coefficient for $x$ achieves same stationary level for both Grid Discretization and Regress Now. Although an exact comparison is  not possible between Regress Now and Regress Later, we continue to observe similar sign and dynamics for each of the coefficients. However, getting away with almost no noise in the dynamics of the estimated coefficients of Regress Later compared to Regress Now is essentially magical.
 
%  Note that for Regress Later the coefficients converge to different values than Regress Now. In order to compare the two methods we should point out that they approximate the same quantity only when basis functions $\hat{\phi}$ are used for the former and $\phi$ for the latter. The conditional expectation $\hat{\phi}$ can be computed as follow:
%  \[
%  \begin{split}
%  \hat{\phi}_k(x,i,d)
% & =\mathbb{E}\big[X^p_{n+1},I^q_{n+1}\big|X_n=x,I_n=i,d_n=d\big]\\
% & =I_{n+1}^q\sigma^p dt^{\frac{p}{2}}\sum_{k=0}^p\mathbb{I}_{\{(p-k)\text{ is odd}\}}{p \choose k}\Big(x\frac{1-\lambda dt}{\sigma\sqrt{dt}}\Big)^k\prod_{j=1}^{\frac{p-k}{2}}(2j-1)
% \end{split}
%  \] 
As a result, we define a stationary policy $d(x,w,m)$ to be used in a longer time horizon than the one employed for its estimation which performance are comparable to the time dependent policy $d(t,x,w,m)$. 

We finally tested the value of both stationary and time dependent policy and found that the performance of the stationary policy is comparable to that of the time dependent policy. 

\subsection{Analysis of the controllers}
\label{controller_analysis}

In this section we compare the control policies estimated by the three algorithms and we try to assess whether one of the approaches is preferable.

\subsubsection{Control maps}

\begin{figure}
        \centering
        \begin{subfigure}[b]{0.33\textwidth}
            \centering
            \includegraphics[width=\textwidth]{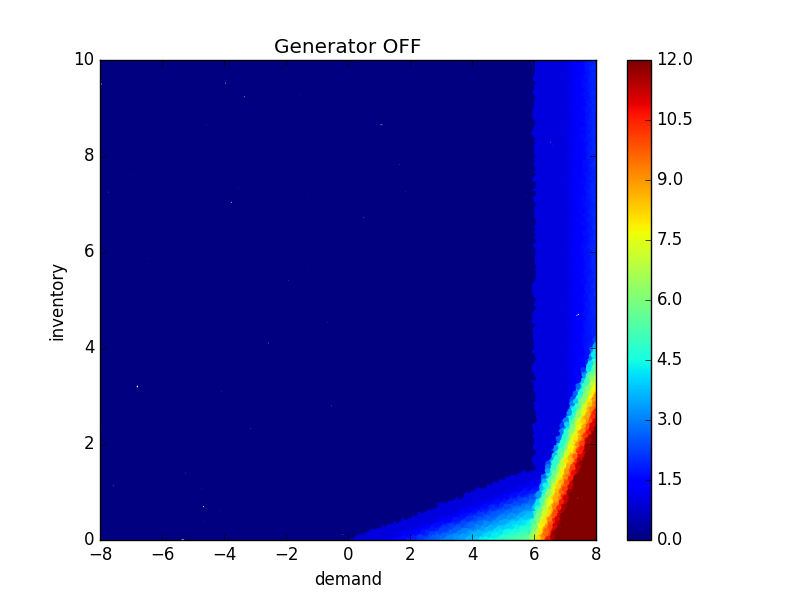}
            \caption[]%
            {{\small }}    
        \end{subfigure}
        \hspace{-40pt}
        \begin{subfigure}[b]{0.33\textwidth}  
            \centering 
            \includegraphics[width=\textwidth]{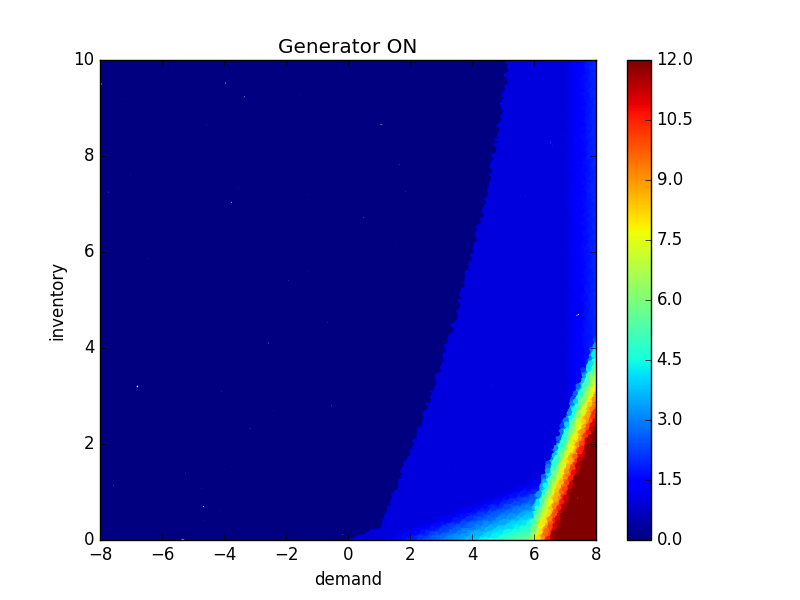}
            \caption[]%
            {{\small }}    
        \end{subfigure}
        \hspace{-30pt}
        \begin{subfigure}[b]{0.33\textwidth}  
            \centering 
            \includegraphics[width=.95\textwidth]{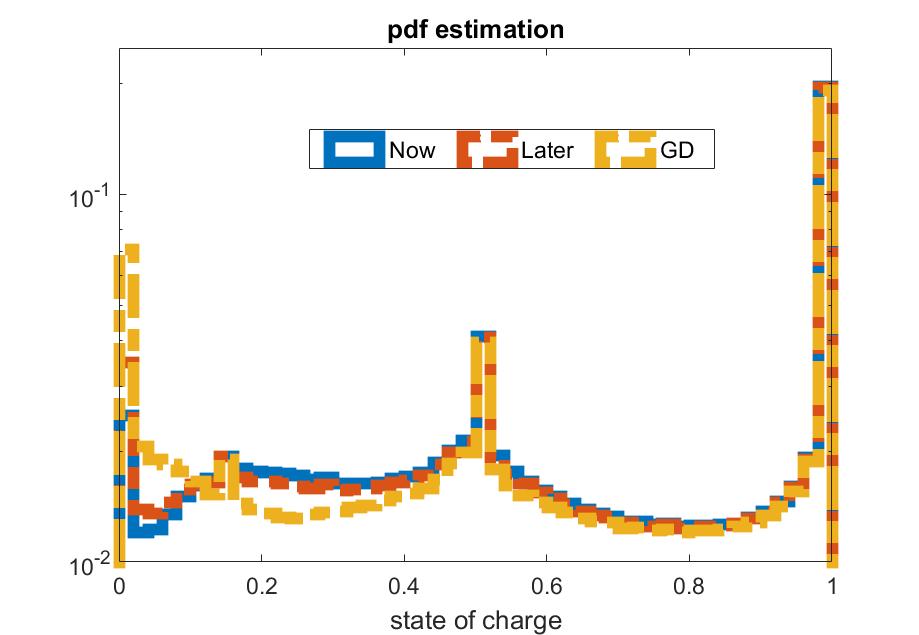}
            \caption[]%
            {{\small }}    
        \end{subfigure}
    \caption{In the figure above we show, in the two left-most panels, an example of control map produced by the Regress Later algorithm. Notice the difference depending on the state of the generator. In the right-most panel we display the estimated probability density function of the state of charge of the battery associated with the use of the three policies. It can be observed that Regress Later and Grid discretization induce very similar distributions.}    
    \label{fig:control_maps_comparison}    
    \end{figure}

We compare now the stationary control policies produced by the different algorithms; recall that these policies are feedback to the state, i.e. can be written as function $d_m(x,w)$. Figure \ref{fig:control_maps_comparison} displays an example of the feedback control policy in the form of control map, a graphical representation of the value of the optimal control for each pair $(x,w)$.

% \begin{figure}
% \hspace{-30pt}\includegraphics[width=0.7\textwidth]{control_maps_temp}
% \hspace{-5pt}\includegraphics[width=0.3\textwidth]{algo_comparison_inv.jpg}
% \caption{ In the figure above we show, in the two left-most panels, an example of control map produced by the Regress Later algorithm. Notice the difference depending on the state of the generator. In the right-most panel we display the estimated probability density function of the state of charge of the battery associated with the use of the three policies. It can be observed that Regress Later and Grid discretization induce very similar distributions. }
% \label{fig:control_maps_comparison}
% \end{figure}

We observed that the three policies agree with the intuition that the diesel generator should produce more power when residual demand is high and inventory is low. We can also notice that the switching cost influences the policy, forcing the diesel to keep running for longer in order to charge the battery sufficiently and avoid turning ON and OFF the generator too often. Just by observation of the control maps little difference can be found among the algorithms, we display in Figure \ref{fig:control_maps_comparison} the effect of the control policy on a the state of charge of the battery. It can be observed from the estimated unconditional probability density of the process $I$ that the policies induced by Regress Now and Regress Later are very similar. Both seem to induce a peculiar mass of probability around $I_n=2.5$, differentiating the behavior of the inventory compared to Grid Discretization. The distribution of the state of charge, obtained by plotting the histogram of all simulations over all time steps, shows that Regress Now and Regress Later does not fully exploit the whole inventory but rather they are more conservative, saving energy to avoid to turn ON the diesel generator in the future. In the next section we will investigate the value associated to this control maps.

\subsubsection{Performance of the policies}

In order to assess the performance of each policy in an unbiased manner, we select a collection of simulated paths of the residual demand process $X$, and record the costs associated with managing the microgrid as indicated by each control map.

\begin{figure}
\includegraphics[width=1\textwidth]{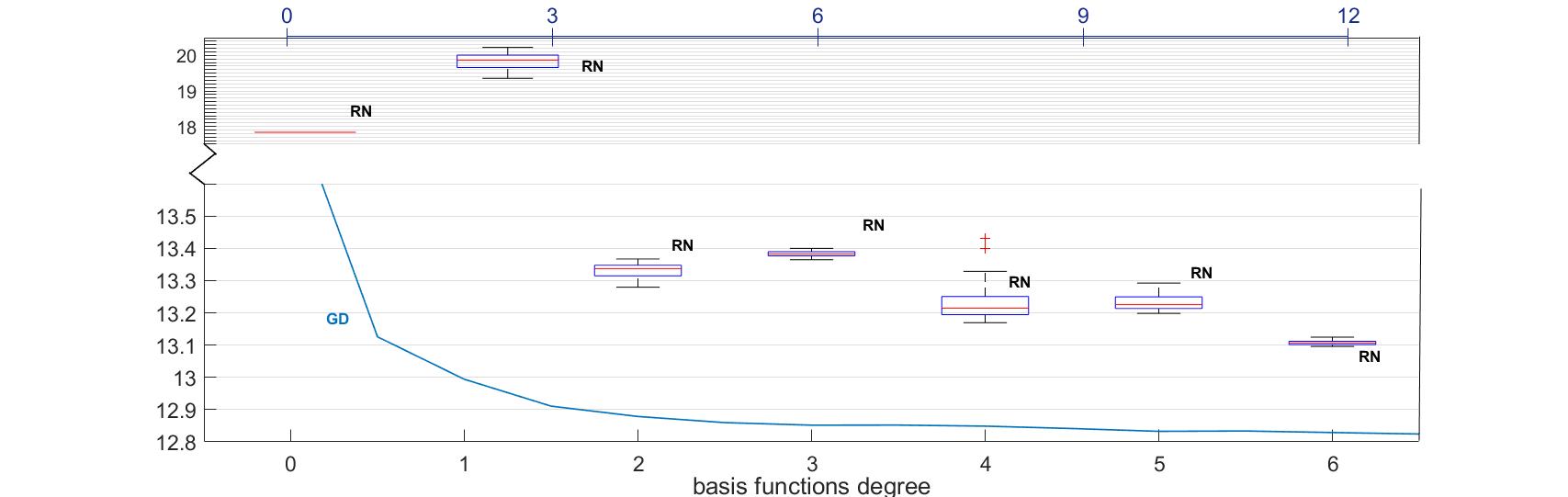}
\caption{The figure in display shows the reduction in operating cost when higher degree polynomials are added to the basis functions, in the case of RN regression, or more inventory points in Grid Discretisation. Notice the peculiar behavior of even/odd degree of basis functions in the RN regressions. Similar analysis was performed for Regress Later and the results are available on request. }
\label{fig:value_comparison}
\end{figure}

We first study how the quality of each policy improves when we increase the computational budget given to each algorithm to compute the stationary policy. In Figure \ref{fig:value_comparison}, we show the estimated value of the policy when the initial state of the system is $(x,i,m)=(0,5,0)$ for polynomial basis functions of increasing degree, for 2D regression. In case of GD we increase the number of discretisation points for the inventory. In particular we make the computational time increase by providing the problem with more training points and more parameters to use in the definition of $\mathcal{C}$ as increasing the number of basis functions. In the case of 2D regression, surprisingly, we noticed that the performance of the estimated control improves only when polynomials of even degree are added, and the effect is more prominent for Regress Later.

We notice from the comparison that Grid Discretisation converges quickly, resulting in the best algorithm in terms of trade off between running time and precision. Among the 2D regressions, we observe similar bias for Regress Now and Regress Later (not displayed in order to maintain clear presentation, but available on request), however latter  has lower standard error. This is not surprising because Regress Later  has only one element of approximation error due to finite basis functions while Regress Now has error attributed to two sources, first, due to finite basis function and second, pathwise estimation of the conditional expectation.  

\subsection{System behavior}
\label{system_behaviour}

In the previous section we selected Grid Discretisation  to be the best performing algorithm by our criteria. In the following we shall always employ Grid Discretisation to conduct our study of the sensitivity of the control policy and the associated cost of managing the grid to some of the parameters of the model. 

The aim of the section is to build a solid understanding of the behavior of the microgrid in order to get an insight into the optimal design of the system. We decided to study the following aspects of the grid: battery capacity, represented by $I_{max}$; different proportion of renewable production, via the volatility $\sigma$ and the mean reversion $b$; tenable behavior of the policy, via the switching cost $K$ and curtailment cost $C$.

In order to be able to carry out our analysis, without introducing cumbersome economic and engineering details regarding the microgrid components, we have to make very simplistic assumptions. Our aim is however to guide the reader through a methodology that can be replicated to study real world microgrid systems.

\subsubsection{Battery capacity}

\begin{figure}
\centering
\includegraphics[width=0.6\textwidth]{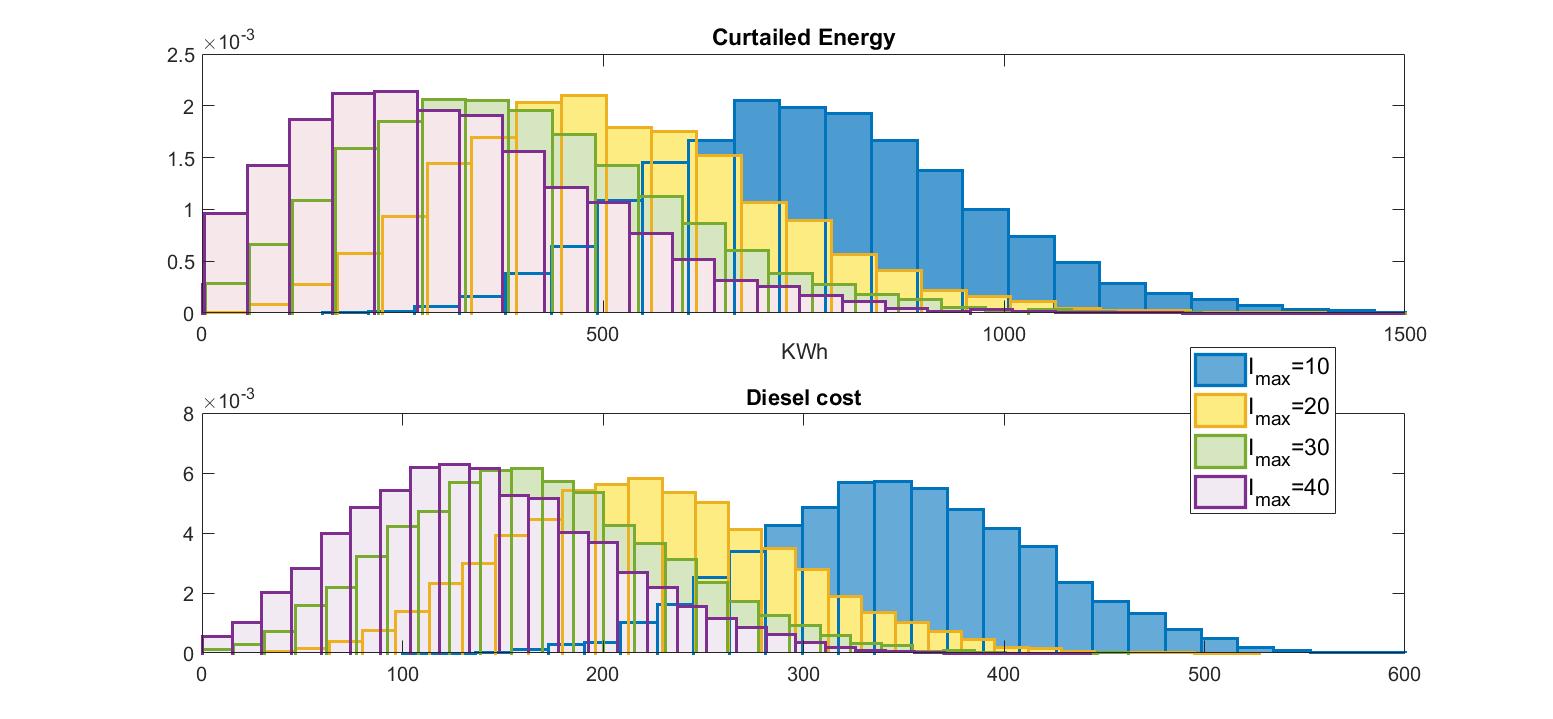}
\caption{In the figure above we show histograms for different levels of battery capacity. In the top panel we display the estimated probability density of the  curtailed energy, while in the bottom panel the estimated density of the cost of operating the diesel generator. Notice that the decrease in cost and curtailed energy per KWh of additional capacity is smaller for high capacity batteries.}\label{fig:sens_battery}
\end{figure}

We study first the behaviour of the system relatively to changes in the capacity of the battery. We would expect to observe negative correlation between the quantity of diesel consumed and the battery size. We display in Figure \ref{fig:sens_battery} both the quantity of energy curtailed and the cost of running the diesel generator for different values of the battery capacity. We can observe that, as expected, increasing the size of the battery leads to lower diesel usage thanks to the higher proportion of renewable energy that is retained within the system. As the capacity of the battery reaches 30/40 KWh, we start observing a decrease in the cost-reduction per KWh of additional capacity suggesting that further analysis should be run in order to understand up to which size it is worth to pay to add storage capacity to the system.

We show now how to infer information about the optimal sizing of the battery, minimizing the trade off between the installation cost of a bigger battery and the reduced use of the diesel generator. Consider however that including battery ageing in the stochastic control problem is outside the scope of this paper but rather in this section we present only a post-optimization analysis.  Assuming that the microgrid runs under similar conditions for the next 10 years, we can quickly estimate the total throughput of energy for the different battery capacities. Consider now that a battery has not an infinite lifetime, but rather it should be scrapped after equivalent 4000 cycles (amount of energy for one full charge and discharge). Under the previous assumptions, we can compute how many batteries would be necessary to cover the next 10 years of operations. Similarly, using the data relative to the usage of diesel generator for different levels of capacity, we can compute the operating cost of the diesel generator over the same time period. Further exploiting the assumption about the lifetime of a battery, we obtain the cost of running the grid for 10 years as a function of the number of batteries. To conclude, assuming a linear cost of 400 $\EUR$/KWh of capacity, we work out the installation cost of the different-size storage devices.

\begin{figure}
\centering
\includegraphics[width=0.8\textwidth]{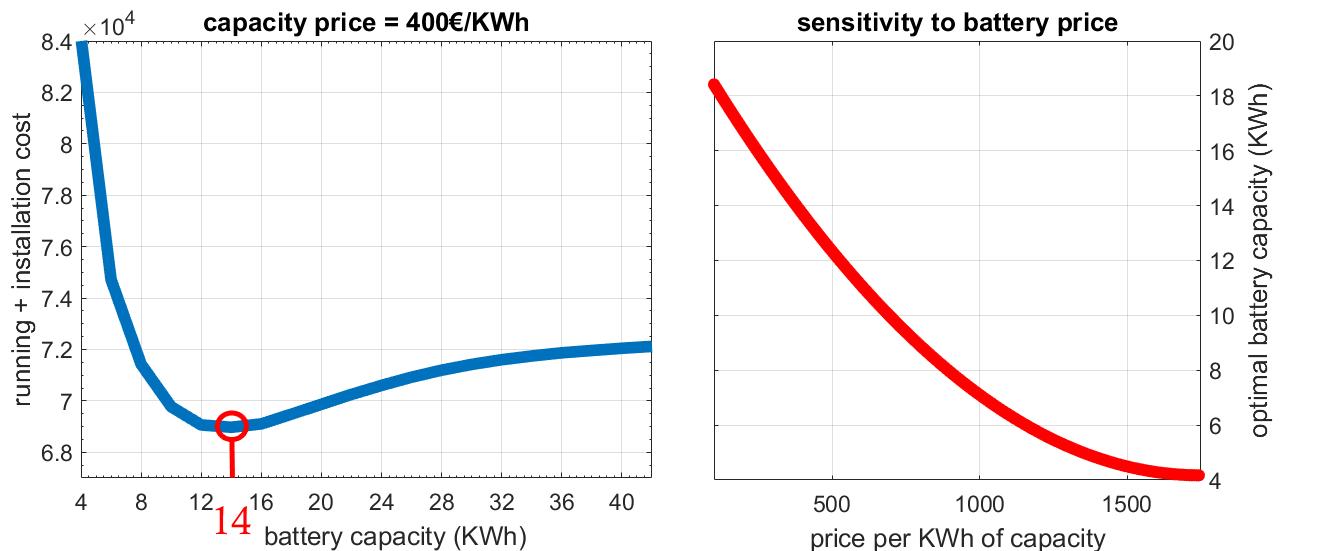}
\caption{In the figure above we compute the total cost of installing and running the grid for ten years, assuming we replace the battery every 4000 cycles, and plot it against the battery capacity (left panel). From the corresponding minimum we can work out the optimal battery capacity and, further, compute the sensitivity of such result with respect to the cost per KWh of capacity.}\label{fig:opt_battery}
\end{figure}

Once this information is collected we search for the minimum of the sum of installation and running cost and, in turn, we compute the optimal capacity. Figure \ref{fig:opt_battery}, on the left, displays a graphical summary of the procedure just described and shows that in our problem the optimal size of the battery is $14$ KWh under the current set of assumptions. Further, we study how much our result is affected by the cost per KWh of capacity, repeating the procedure above. We find that, as expected, as cost increase the size of the optimal battery decreases. Figure \ref{fig:opt_battery}, on the right, displays such behaviour.   

\subsubsection{Renewable penetration}
In this section we want to investigate how robust the microgrid is  to higher penetration of renewable generation, or, in other words, to what extent the algorithm can cope with increasing randomness and decreasing predictability of the system. To model this phenomena we assume that greater penetration of renewables can be modeled by increasing both the parameters for volatility $\sigma$ and the mean reversion rate $\lambda$. Increasing these two parameters makes the problem more difficult to solve, given that the control policy can rely less and less on the statistical properties of the process $X$.  
\begin{figure}
\centering
\includegraphics[width=0.45\textwidth]{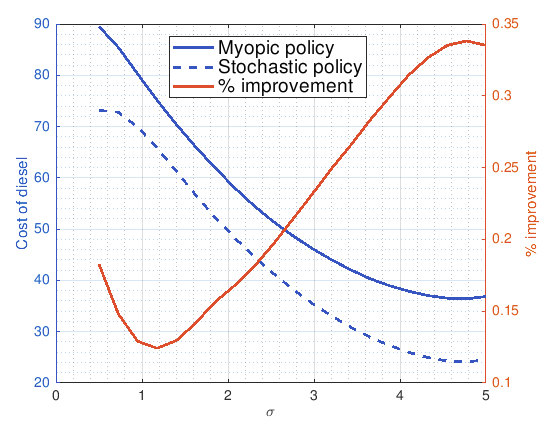}
\caption{The figure represents the cost of the diesel usage for stochastic and myopic policy as a function of $\sigma$. The orange curve represents the percentage improvement in cost due to as a proportion of cost of myopic policy.}\label{fig:renewablePenetration_1}
\end{figure}

%\begin{figure}
%\centering
%\includegraphics[width=0.7\textwidth]{robust}
%\caption{\ab{I would remove it, indeed we have no replacement as today} In the picture above, in the bottom panel, we show a set of simulated trajectories for increasing volatility and mean reversion. It can be observed that both roughness and absolute value increase with increase in renewable penetration. In the top panel we show the relative increase in cost as function of the increase in cumulative demand. Notice that the very slow increase, partially explained by the corresponding increase in cumulative energy produced, shows that our algorithm can deal with high levels of renewable penetration.  }\label{fig:sens_ren}
%\end{figure}

    \begin{figure}
        \centering
        \begin{subfigure}[b]{0.45\textwidth}
            \centering
            \includegraphics[width=\textwidth]{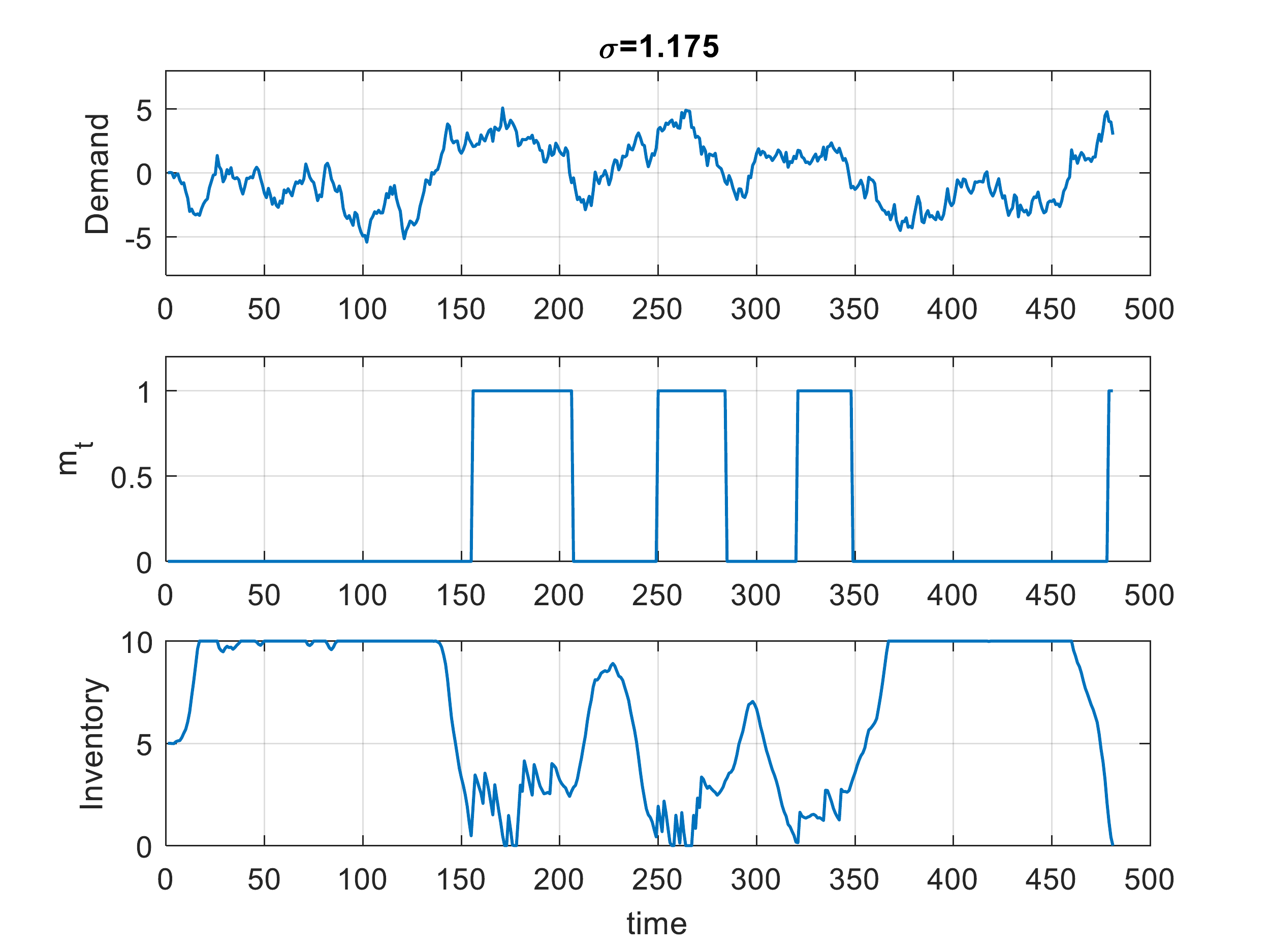}
            \caption[]%
            {{\small }}    
        \end{subfigure}
        \quad
        \begin{subfigure}[b]{0.45\textwidth}  
            \centering 
            \includegraphics[width=\textwidth]{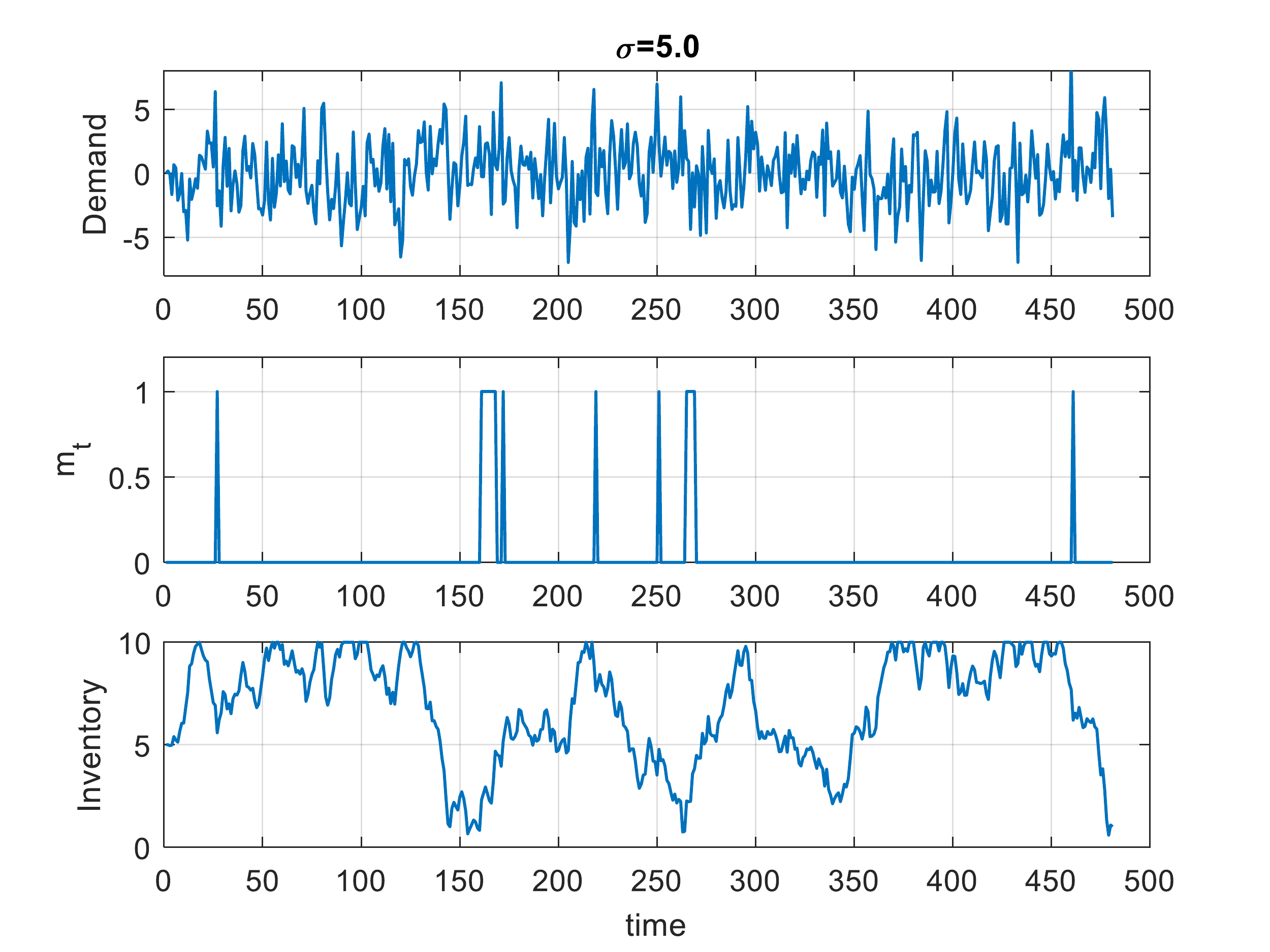}
            \caption[]%
            {{\small}}    
        \end{subfigure}
    \caption{Figure in the left and right panel represents demand, diesel usage and the inventory dynamics for low and high $\sigma$ respectively. It is important to mention that the mean reversion rate was chosen as $\lambda := \sigma^2/8$, in order to ensure a constant volatility of the process regardless of $\sigma$. Notice the low usage of the diesel generator in the figure on the right compared to the one on the left.}
    \label{fig:renewablePenetration_2}    
    \end{figure} 

In order to establish the real added value provided  by our stochastic optimization algorithm, we compare the estimated policy with an heuristic myopic control which can be reproduced in our model solving the dynamic programming equation \eqref{eq:dpe} taking constant conditional expectation with respect to the control. We show the value of the two control policies as function of the increasing learning difficulty in Figure  \ref{fig:renewablePenetration_1} where we observe that the value of accounting for statistical estimation of future conditional expectations when taking decisions decreases.

In figure \ref{fig:renewablePenetration_1} we present cost of diesel as a function of $\sigma$ for stochastic and myopic policy. Since increasing $\sigma$ alters the volatility of the distribution, we define the mean reversion rate $\lambda := \sigma^2/(2c)$ in order to ensure that the volatility of the process is constant while we increase $\sigma$. The stochastic policy leads to at least $12\%$ reduction in the cost of the diesel usage, compared to the myopic policy, and the difference magnifies with increasing ``fluctuations" in the process. The decreasing relationship of the cost with $\sigma$ signifies the importance of the battery storage system in the microgrid which absorbs the sharp change in the demand. In figure \ref{fig:renewablePenetration_2} we compare the demand for two different levels of the $\sigma$, the dynamics of the diesel generator and the inventory. Notice significantly less usage of the diesel for high fluctuations, $\sigma=5$, compared to $\sigma=1.175$. 

The results of this experiment are affected by the over-pessimistic assumption of modeling greater penetration of renewables with an increasingly unpredictable, and eventually completely random, residual demand process. This sort of analysis can however provide insight into how much (weather and load) forecasting capability will be necessary for a given level of renewable penetration.

\subsubsection{Switching and curtailment}
\begin{figure}
\centering
\includegraphics[scale=0.3]{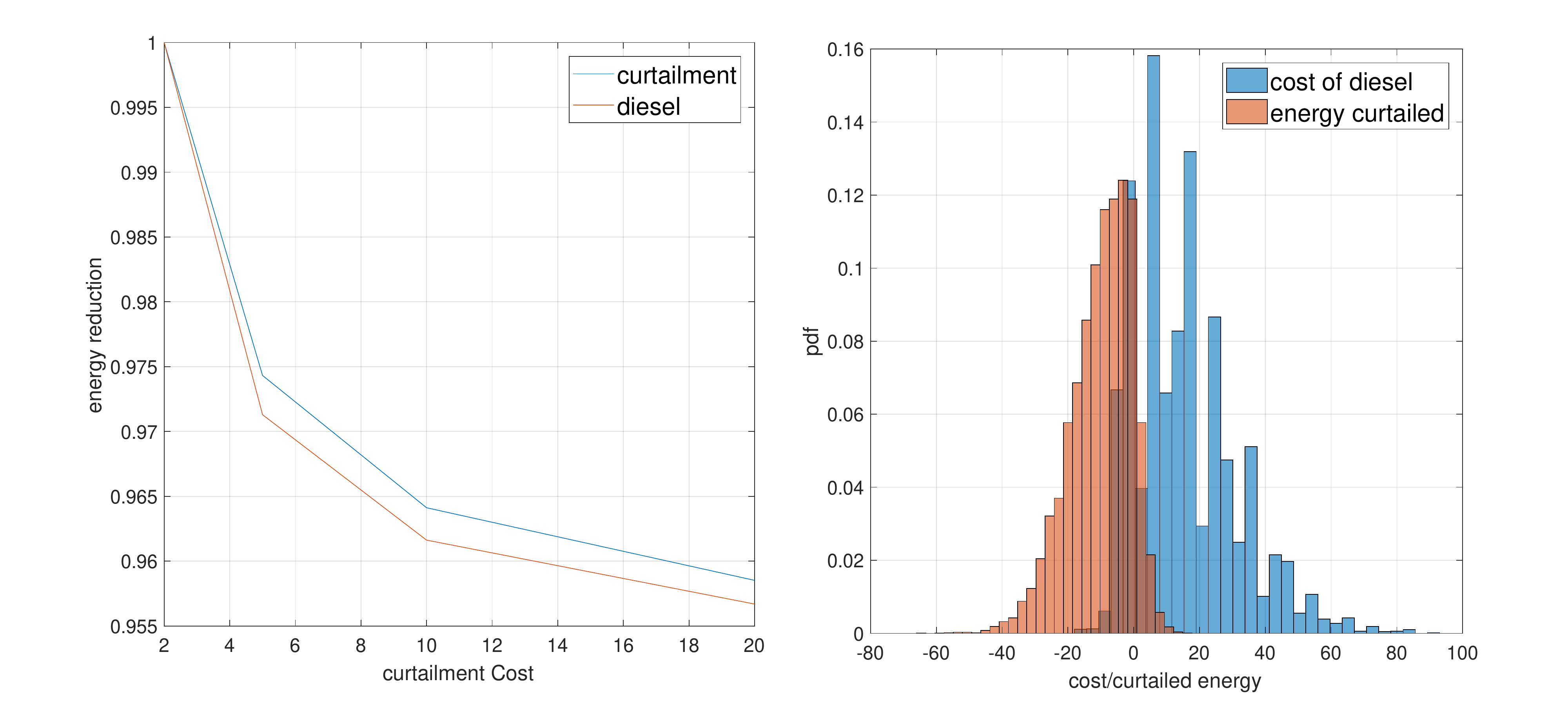}
\caption{Line plot on the left, represents the impact of curtailment cost on the total curtailed energy for different C as a proportion of curtailed energy at C=2. The histogram  on the right, represents the difference in cost of diesel and the curtailed energy for C=20 and C=2. Notice the increase in curtailment cost leads to reduced curtailed energy but at the expense of inefficient diesel usage.}
\label{fig:curtailment}
\end{figure}

We conclude this section by analyzing the dependence of the system behavior on two key parameters in the model: switching cost $K$ and curtailment cost $C$. Switching cost is a system's property and the microgrid controller has little freedom over, however the controller can significantly reduce the amount of curtailed energy by choosing the appropriate curtailment cost. In figure \ref{fig:curtailment}, we observe that increasing the curtailment cost reduces the total curtailed energy by approximately 4\%. However, it comes at the cost of inefficient usage of the diesel generator, which is represented on the right in the figure \ref{fig:curtailment}. The histograms represent the difference between the cost of diesel usage (blue) and the energy curtailed (orange) for C=20 and C=2. Positive diesel cost depicts inefficient usage of the diesel at C=20 compared to C=2. Depending upon the specific cost functional for the diesel, the controller can use C as a parameter for better optimization.

%    \begin{figure}
%        \centering
%        \begin{subfigure}[b]{0.4\textwidth}
%            \centering
%            \includegraphics[width=\textwidth]{curtailedEnergy.pdf}
%            \caption[]%
%            {{\small }}    
%        \end{subfigure}
%        \quad
%        \begin{subfigure}[b]{0.4\textwidth}  
%            \centering 
%            \includegraphics[width=\textwidth]{costFigure.pdf}
%            \caption[]%
%            {{\small }}    
%        \end{subfigure}
%    \caption{Line plot on the left, represents the impact of curtailment cost on the total curtailed energy for different C as a proportion of curtailed energy at C=2. The histogram  on the right, represents the difference in cost of diesel and the curtailed energy for C=20 and C=2. Notice the increase in curtailment cost leads to reduced curtailed energy but at the expense of inefficient diesel usage.}
%    \label{fig:curtailment}    
%    \end{figure} 

\begin{figure}
        \centering
        \begin{subfigure}[b]{0.4\textwidth}
            \centering
            \includegraphics[width=\textwidth]{RL_q1_controlMap.png}
            \caption[]%
            {{\small $K=2$}}    
        \end{subfigure}
        \quad
        \begin{subfigure}[b]{0.4\textwidth}  
            \centering 
            \includegraphics[width=\textwidth]{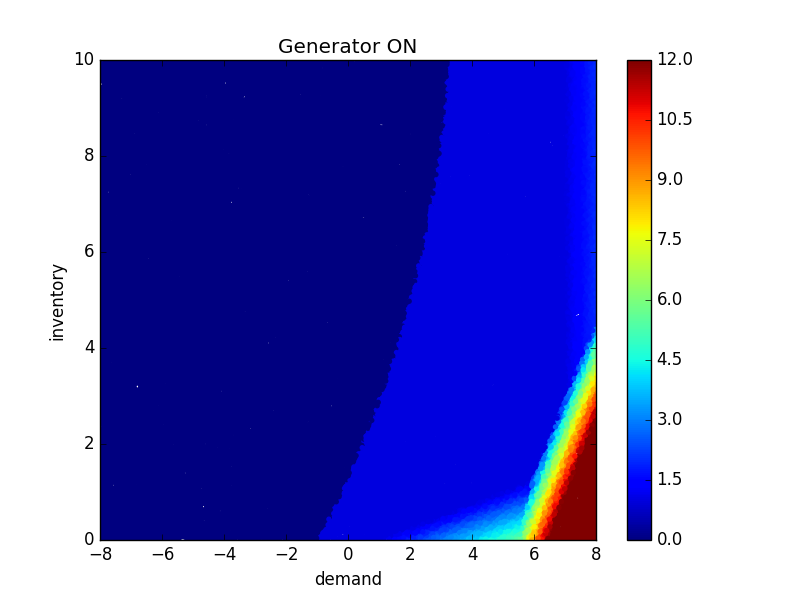}
            \caption[]%
            {{\small $K=5$}}    
        \end{subfigure}
    \caption{Figure on the left represents the control map for switching cost $K=2$, while the figure on the right represents the control map for $K=5$ when the generator is ON. Notice the increase in area for light blue (corresponding to $d=1$) in the figure on the right because of increased switching cost. }
    \label{fig:switchingCost}    
    \end{figure}

The optimal policy when the generator is ON $m_t=1$ is significantly altered depending upon the switching cost. For example, in figure \ref{fig:switchingCost}, we present the control maps associated with K=2 and K=5. As expected, larger switching cost disincentivise  the controller to switch OFF the diesel generator once it's ON. However, we don't observe "significant" change in the control policy due to increase in switching cost when the generator is OFF.

\section{Comparison with deterministically trained policy}
\label{deterministic_comparison}
In this section we compare our stochastic optimization algorithm with a deterministically trained policy. The latter is widely used in online optimization where the solution is computed with respect to the best forecast available at a given time. We emulate this situation by computing the optimal set of actions for a particular deterministic demand trajectory at different levels of the inventory. We assume that the forecast of the demand is given by:
\begin{equation}
\label{eq:deterministic_demand_dynamics}
X_{t+1}=  X_{t} + 0.5(6\sin(\frac{\pi t}{12}) -X_{t} )\Delta t  ; \quad t\in\{0,1,\dots,T-1\}.
\end{equation}
Equation \eqref{eq:deterministic_demand_dynamics} implies periodicity of one day in the residual demand and is equivalent to $\sigma = 0$, $b=0.5$ and $\Lambda_t = 6\sin(\frac{\pi t}{12})$1 in \eqref{eq:residual_demand_dynamics}. Zero volatility in the residual demand curve leads to a deterministic optimal control problem, rather than a stochastic control problem we have presented in section 5. 

Notice that the deterministic optimal control problem results in a sequence of control maps $d_t:(w,m) \to\mathbb [d_{min},d_{max}]\cup 0$. As a result, although the policy has been trained on a deterministic residual demand,  it dynamically adapts itself to different inventory levels and state of the diesel generator, when tested in a stochastic environment. We present the modified algorithm in  \ref{algo:deterministic}. There are two key differences from the previous algorithm, first, we use one dimensional projection of the value function and second, we replace regression with interpolation since there is no randomness left in the problem. 

\begin{algorithm}
\caption{Regression Monte Carlo algorithm for deterministic demand}
% \textbf{input:} \ab{operating parameters?} forecast $\Lambda$, .
\begin{algorithmic}[1]

\State Simulate $\{X_{t}\}_{t=1}^{N}$ according to its dynamics;
\State  Discretize $I_t$ into $M$ levels indexed by $j$ s.t. $\{I_{t}^{j}\}_{j=1}^{M}$ ;

\State Initialize the value function $V(T, I_{T}^j,m_T) =g(I_T^j), \quad \forall j=1,\,\dots,\, M$ and $m_T=\{0,1\}$ ;
\For{$t=N-1$ to $1$}
    \State Find interpolation function  $\mathcal{B}(t+1,I_{t+1},m)$ for $\{ V(t+1, I_{t+1}^j,m_{t+1}) \}_{j=1}^{M} \quad \text{for each } m=0,1$ 
    \State Compute the set of admissible controls as $\mathcal{U}_{t}$
    \For{$j=1$ to $M$}
    \For{$m=0$ to $1$}
\State $F = \mathcal{B}(t+1,I_{t}^j,0)$
\State 
\[
V(t,I_{t}^j,m)=
\begin{cases}
\min\limits_{d\in\mathcal{U}_{t} \setminus \{0 \} }\Big\{ p\rho(d)+ CS_t\mathds{1}_{\{S_t<0\}} + \mathcal{B}(t+1,I_{t}^j-B_t^d,1) \Big\}+K\mathds{1}_{\{m=0\}} \wedge F  &\quad \text{if } 0 \in \mathcal{U}_{t}\\ 
\min\limits_{d\in\mathcal{U}_{t} }\Big\{ p\rho(d)+ CS_t\mathds{1}_{\{S_t<0\}} + \mathcal{B}(t+1,I_{t}^j-B_t^d,1) \Big\}+K\mathds{1}_{\{m=0\}} &\quad \text{otherwise}
\end{cases}
\]
\EndFor
    \EndFor
\EndFor
\end{algorithmic}
\textbf{output:} control policy $\{\mathcal{B}{(t,\cdot,\cdot)}\}_{t=2}^{N}$.
\label{algo:deterministic}
\end{algorithm}

In order to understand the solution of the deterministic problem, in figure 13 we present the dynamics of the optimal control and inventory corresponding to the demand faced in (A). As expected, diesel switches on when the demand is high and it keeps it running just long enough that the battery is empty before it faces negative residual demand to charge the battery. Moreover, there is substantial curtailment of energy since the battery is not large enough to store all the  excess energy.

    \begin{figure}
        \centering
        \begin{subfigure}[b]{0.3\textwidth}
            \centering
            \includegraphics[width=\textwidth]{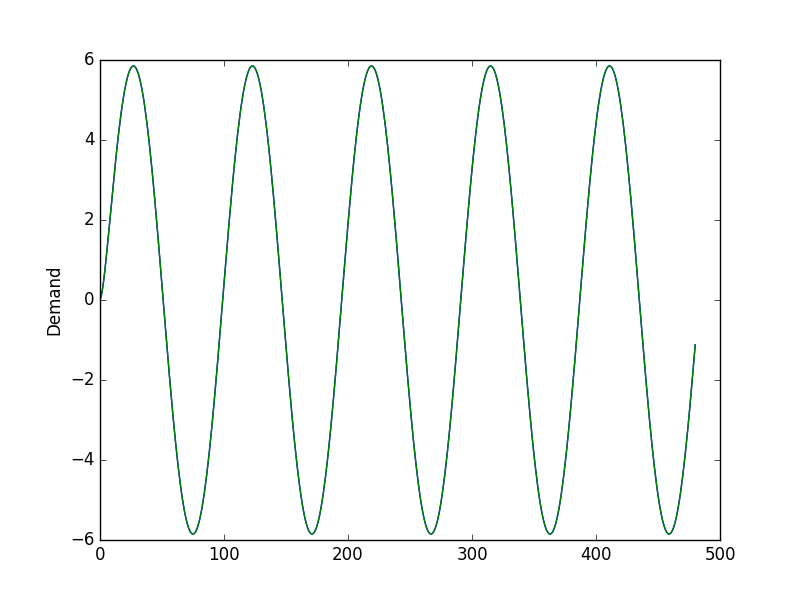}
            \caption[]%
            {{\small Demand}}    
        \end{subfigure}
        \hfill
        \begin{subfigure}[b]{0.3\textwidth}  
            \centering 
            \includegraphics[width=\textwidth]{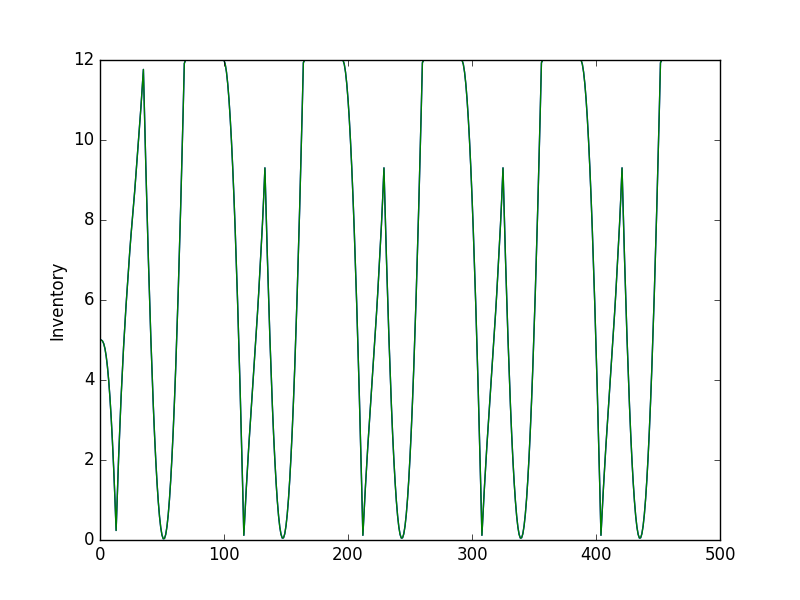}
            \caption[]%
            {{\small Inventory}}    
        \end{subfigure}
        \hfill
        \begin{subfigure}[b]{0.3\textwidth}  
            \centering 
            \includegraphics[width=\textwidth]{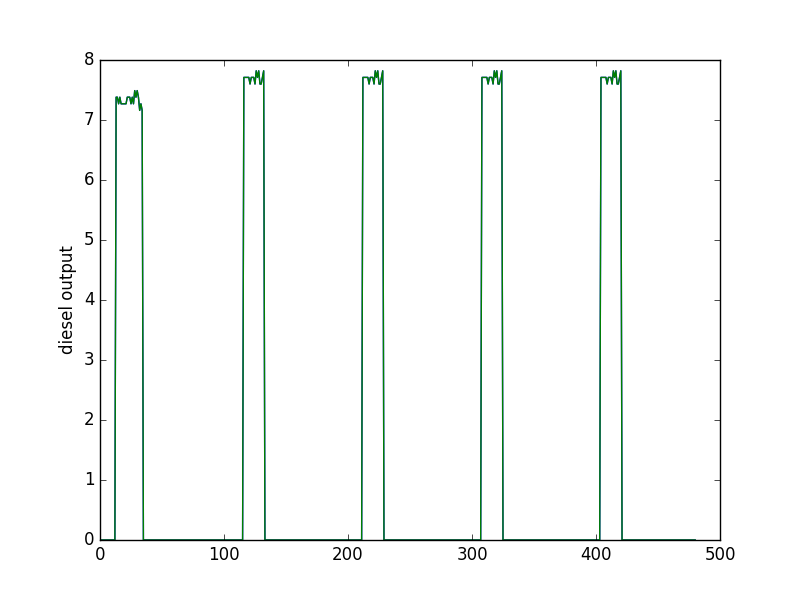}
            \caption[]%
            {{\small Diesel Output}}    
        \end{subfigure}
	        \caption{The image illustrates the dynamics of the inventory and control  for the deterministic control problem. Figure (A) represents the demand in equation \eqref{eq:deterministic_demand_dynamics}, the optimal control of the diesel in figure (C) and the corresponding dynamics of the inventory in figure (B).}
    \end{figure}

In order to quantify the gain due to formulating the microgrid management problem as a stochastic control rather than traditional deterministic control, we compare the performance of the deterministically trained strategy of this section to its stochastic counterpart developed in this paper. While the deterministic control problem was solved using the residual demand curve \eqref{eq:deterministic_demand_dynamics}, the stochastic control problem was fed in with the residual demand curve \eqref{eq:residual_Sinusoidal_demand_dynamics}. Finally, we test both the strategies on fresh out-of-sample paths following the residual demand  \eqref{eq:residual_Sinusoidal_demand_dynamics}.
\begin{equation}
\label{eq:residual_Sinusoidal_demand_dynamics}
X_{t+1}=  \Big(X_{t}+0.5(6\sin(\frac{\pi t}{12}) -X_{t} )\Delta t + 2 \sqrt{\Delta t}\xi _{t}\Big)\;  \wedge 10\; ; \quad t\in\{0,1,\dots,T-1\}
\end{equation}

\begin{figure}
\begin{floatrow}
\ffigbox{%
  \includegraphics[scale=0.25]{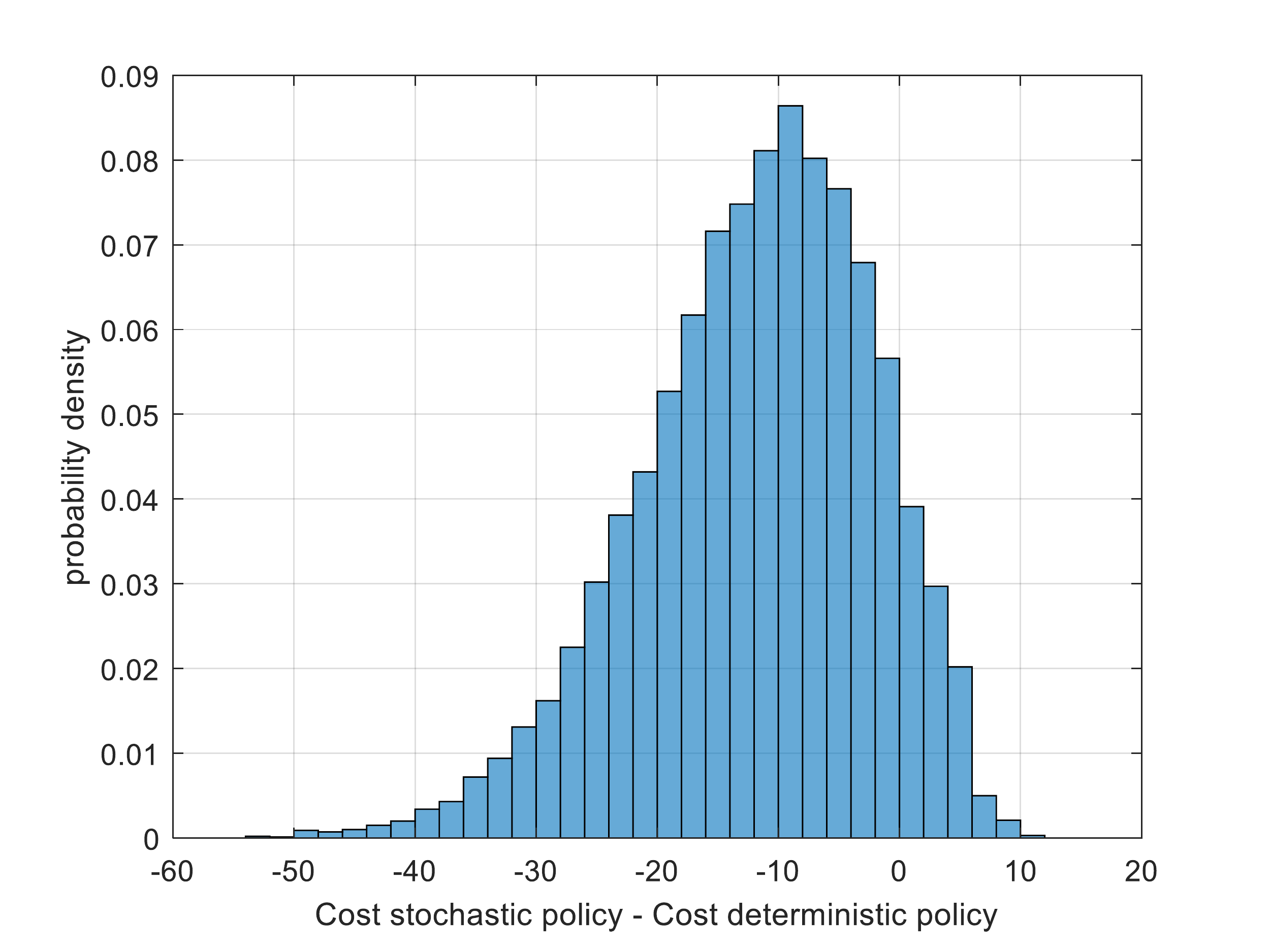}%
}{%
  \caption{Difference of the Cost of Stochastic and deterministic policy for K=5}
     \label{fig:histogramCost}
}
\capbtabbox{%
\begin{tabular}{|l|r|r|r|}
\hline
Switching Cost & K=2 & K=5 & K=10 \\ \hline
Deterministic & 138.56 & 162.63 & 201.52 \\ \hline
Stochastic & 131.86 & 150.49 & 178.22 \\ \hline
\% difference & 4.84\% & 7.46\% & 11.56\% \\ \hline
\end{tabular}
}{%
  \caption{Comparison of deterministic and stochastic trained policy.}
   \label{table:compare_DeterministicStochastic}
}
\end{floatrow}
\end{figure}

In figure \ref{fig:histogramCost}, we present the histogram of the cost from the stochastic policy and the deterministic policy pathwise for 10,000 out-of-sample paths. As evident, most of the distribution lies on the negative side, implying gain due to stochastic policy. To measure this difference, in table \ref{table:compare_DeterministicStochastic}, we quantify the gain of the stochastic policy for different switching cost. For switching cost of K=5, we observe that the stochastic policy is 7.5\% better than the deterministic policy. As the switching cost increases, mistakes made by deterministic policy become more expensive leading to higher percentage difference. 
%\iffalse
%\fi

% use the same trajectory as mean reversion level $\Lambda_t$, introduced in \eqref{XX}, for the stochastic dynamics of the residual demand process $X$ used by our stochastic algorithm. 

% In order to adapt a deterministic sequence of actions to control policies (functions) needed in the bigger state space in the stochastic framework, we need to modify our estimation of the deterministically trained control as follow:
% \ab{discuss, how do we adapt our decision to different values of X (we interpolate only in I but do not adapt the control to X)?}

% To compare the two approaches we run forward stochastic simulations for the following set of parameters:
% \ab{parameters}

\begin{figure}
\centering
\includegraphics[scale=0.4]{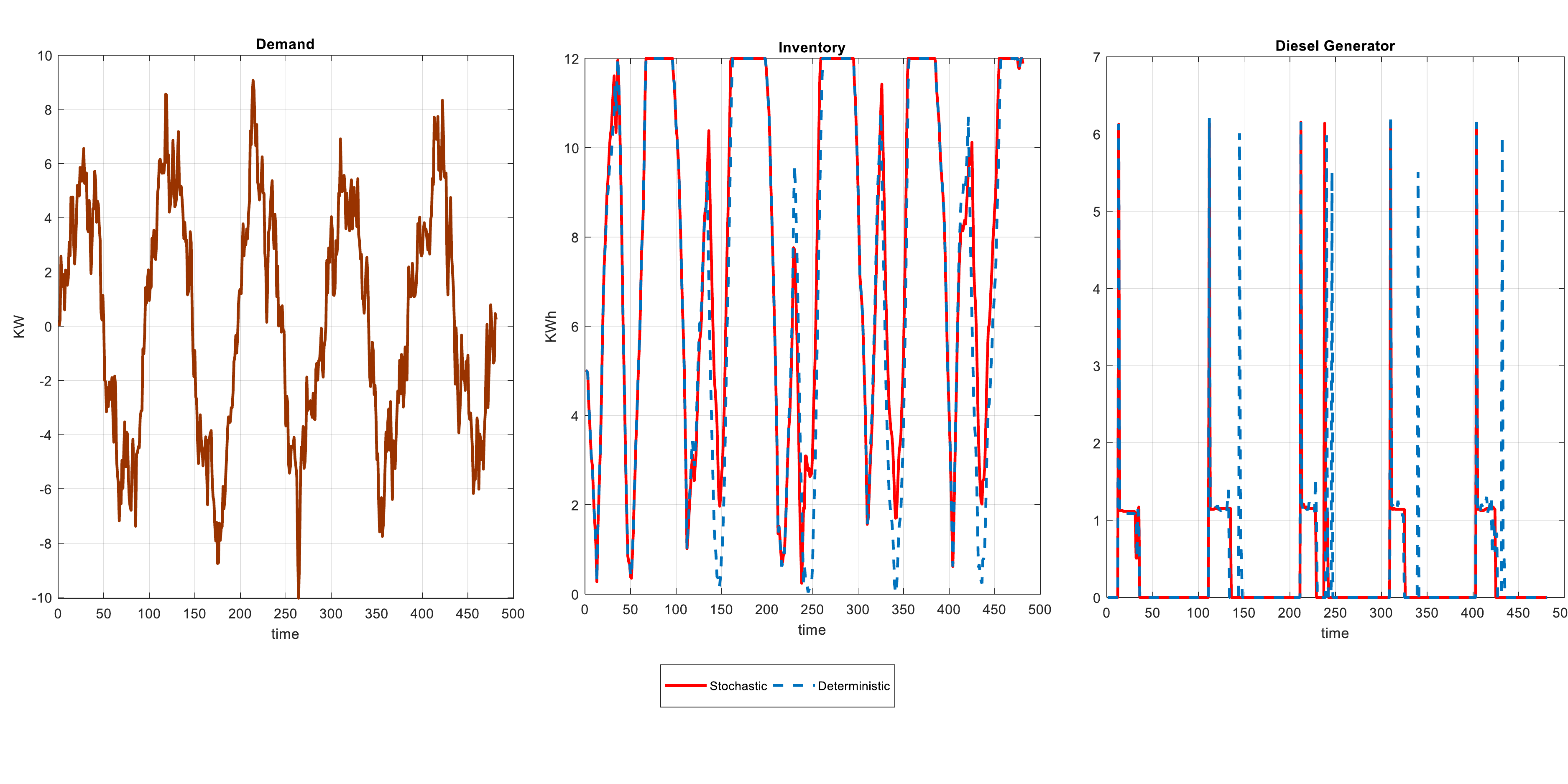}
\caption{The figure above presents the pathwise comparison of stochastic and deterministic policy for the same demand on the left panel. The center panel represents dynamics of the inventory due to control on the right panel. Particularly notice the difference in switching times for the diesel in the deterministic policy and stochastic policy. }\label{fig:deterministic}
\end{figure}

Finally, Figure \ref{fig:deterministic} displays the behavior of inventory and the cost along a random trajectory of residual demand. In blue we show the stochastically trained control policy and in orange the deterministically trained. The stochastic policy has lesser switch of the diesel generator and thus lower costs. The spikes in the cost function for the deterministic policy is due to poor management of the inventory  and thus inefficient usage of the microgrid.

%\ab{discuss conclusion on the comparison}

\section{Conclusion}
In this paper we solved the problem of optimal management of a
microgrid by employing three algorithms from the Regression Monte
Carlo literature, namely: Regress Now, Regress Later and Inventory
Discretization. We find that Inventory Discretization significantly
outperforms the other two methods. Besides algorithm design, we
propose a methodology to optimize the design of the grid and determine the optimal sizing of the battery. In addition, we perform a thorough sensitivity analysis to some of the key parameters, showing the robustness of our solution.  Finally, we compare the control policy estimated by our algorithm to industry standard deterministic control, observing a 5-10\% reduction in cost.

Future research in this direction will include further studies of the
optimal sizing of the battery by explicitly incorporating the wearing
off caused by usage. Another more challenging direction is to
understand the impact of delay, e.g., in the switching of the diesel
generator, on the optimal management of the microgrid. This problem
introduces several mathematical and algorithmic issues which are
currently the focus of our research. 

\section{Acknowledgements}
This research was supported by the FIME Research Initiative. The research of C. Alasseur and X. Warin has also benefited from support by the ANR project CAESARS  (ANR-15-CE05- 0024). The research of Peter Tankov has also benefited from support by the ANR project FOREWER  (ANR-14-CE05- 0028). The research of Aditya Maheshwari has benefited from support by the grant AMPS-1736439. 

\clearpage

\bibliographystyle{apalike}
\bibliography{references}

\end{document}